\documentclass{amsart}
\usepackage{amsfonts}
\usepackage{amssymb}
\usepackage{amssymb}
\usepackage{enumerate}
\usepackage{latexsym}
\usepackage{diagrams}
\newarrow{Mapsto}{|}{-}{-}{-}{>}
\usepackage{nomencl,mathrsfs}

\newtheorem{theorem}{Theorem}[section]
\newtheorem{lemma}[theorem]{Lemma}
\newtheorem{proposition}[theorem]{Proposition}

\newtheorem{corollary}[theorem]{Corollary}

\theoremstyle{definition}
\newtheorem{definition}[theorem]{Definition}
\newtheorem{remark}[theorem]{Remark}

\numberwithin{equation}{section}

\newcommand{\CP}{\mathbb{CP}}
\newcommand{\N}{\mathbb{N}}

\newcommand{\R}{\mathbb{R}}
\newcommand{\Sp}{\mathbb{S}}
\newcommand{\T}{\mathbb{T}}

\newcommand{\co}{\mskip0.5mu\colon\thinspace}
\newcommand{\bas}{{\textup{\scriptsize bas}}}
\newcommand{\ibas}{{\textup{\scriptsize ibas}}}
\newcommand{\inv}{{\textup{\scriptsize inv}}}
\newcommand{\inertia}[1]{{\Lambda}#1}
\newcommand{\cartan}{\mathsf{T}}
\newcommand{\SO}[1]{\mathsf{SO}(#1)}
\newcommand{\Specr}{\operatorname{Spec}_\textup{r}}
\newcommand{\calC}{\mathcal{C}}
\newcommand{\sfB}{\mathsf{B}}
\newcommand{\sfG}{\mathsf{G}}
\newcommand{\id}{\operatorname{id}}
\newcommand{\dec}{\mathscr{D}}

\begin{document}

\title{Stratifications of Inertia Spaces of Compact Lie Group Actions}

\author{Carla Farsi}
\address{Department of Mathematics, University of Colorado at Boulder, Campus
Box 395, Boulder, CO 80309-0395 }
\email{farsi@euclid.colorado.edu}

\author{Markus J. Pflaum}
\address{Department of Mathematics, University of Colorado at Boulder, Campus
Box 395, Boulder, CO 80309-0395 }
\email{pflaum@Colorado.EDU }

\author{Christopher Seaton}
\address{Department of Mathematics and Computer Science,
Rhodes College, 2000 N. Parkway, Memphis, TN 38112}
\email{seatonc@rhodes.edu}

\dedicatory{Dedicated to David Trotman on the occasion of his 60th birthday.}



\begin{abstract}
We study the topology of the inertia space of a smooth $G$-manifold
$M$ where $G$ is a compact Lie group.  We construct an explicit
Whitney stratification of the inertia space, demonstrating that the
inertia space is a triangulable differentiable stratified space. In
addition, we demonstrate a de Rham theorem for differential forms
defined on the inertia space with respect to this stratification.
\end{abstract}

\maketitle

\tableofcontents

\section{Introduction}
\label{sec-Intro}

Let $G$ be a compact Lie group acting (from the left) on a smooth manifold $M$.
In the case where $G$ acts locally freely on $M$, the orbit space $X:=G\backslash M$
is an orbifold. Moreover, in this situation, the inertia space $\inertia{X}$
of the orbifold $X$ can be defined as the quotient of the disjoint union
$\bigsqcup_{g\in G} M^g$ of the fixed point manifolds $M^g$ by the natural action of the
Lie group $G$. It turns out that $\inertia{X}$ is an orbifold as well
which in general has several connected components of varying dimension.
The inertia space of an orbifold has originally been introduced by
{Kawasaki} in \cite[p.~77]{kawasaki1} and subsequently used in
\cite{kawasaki2, kawasaki3}. In these papers, the inertia orbifold served as
a bookkeeping device for the formulation of the topological index in
an orbifold signature theorem resp.~an orbifold index formula.
Since {Kawasaki}'s work, the inertia orbifold has played a major role for all
formulations of index theorems on orbifolds;
see e.g.~\cite{FarsiKIndex,FarsiRelIndex, PflaumAlgIndex,VergneEquivarIndex}.
In addition, the inertia orbifold has been widely studied in connection with the Chern
character for orbifolds, which provides an isomorphism between the (rationalized) orbifold
$K$-theory and the cohomology of the inertia orbifold, as well as with Chen--Ruan orbifold
cohomology, which is additively isomorphic to the cohomology of the inertia orbifold; see
e.g.~\cite{ademleidaruan,BaumConnes,baumbryMacP}.

In the general case, where the action of $G$ is no longer assumed to be locally free,
the space $G\backslash M$ is not necessarily an orbifold but rather a differentiable
stratified space; see \cite[Chap.~4]{PflaumBook}.  In this case, however, an analog
of the inertia orbifold has appeared in connection with the study of the convolution
algebra $\mathcal{C}^\infty (M) \rtimes G$ from the point of view of noncommutative geometry.
More precisely, {Brylinsky} demonstrated in \cite{bryl87}
that the Hochschild cohomology of the convolution algebra
$\mathcal{C}^\infty (M) \rtimes G$ is isomorphic to the space
of relative basic differential forms on an appropriately defined space which
in this paper we will identify with the inertia space of the groupoid
$G\ltimes M$.
Similarly, {Block} and {Getzler} proved in \cite{BlockGetz} that
the periodic cyclic cohomology of the convolution algebra $G $ is isomorphic to a
sheaf of equivariant differential forms defined as a sheaf on $G$.  When $G$ is connected,
this sheaf can also be understood as a sheaf of equivariant relative differential
forms on the inertia space.
In their paper \cite{LupUri}, {Lupercio--Uribe} defined for every
topological groupoid $\mathsf{G}$ an inertia groupoid
$\inertia{\mathsf{G}}$. More precisely, the inertia groupoid $\inertia{\mathsf{G}}$
is  the transformation groupoid $\mathsf{G} \ltimes \mathsf{B}_0$, where
$\mathsf{B}_0$ is the space of loops of the groupoid $\mathsf{G}$, i.e.~the set of all arrows $g$ such
that the source $s(g)$ coincides with the target $t(g)$. The inertia space of $\mathsf{G}$ then
is orbit space of the inertia groupoid  $\inertia{\mathsf{G}}$ or in other words
the quotient space $\mathsf{G}\backslash \mathsf{B}_0$.
If $\mathsf{G}$ is a proper \'etale Lie groupoid representing
an orbifold $X$, the thus defined inertia space coincides with  the inertia orbifold of $X$ as
defined originally by {Kawasaki} and subsequent authors.
See also \cite{ademgomez} for recent results on the $K$-theory of inertia spaces
of  compact Lie group actions where the fundamental group of the Lie group is
torsion-free and all isotropy groups have maximal rank.

With this paper, we aim at defining a general notion of the inertia space of
a proper Lie groupoid and studying its fundamental properties in the basic situation
where the Lie groupoid is a transformation  groupoid  $G\ltimes M$ with $G$ a compact Lie group.
Under this hypothesis, we give an explicit stratification of the inertia space in
Theorem \ref{thrm-StratConnected}.  Additionally, we demonstrate a de Rham theorem
for differential forms on the inertia space in Theorem \ref{thrm-deRhamInertia}.
Note that (locally), the inertia space is a subanalytic set, hence is known to admit a
stratification by \cite{MatumotoShiota}.  However, the stratification
constructed here is given explicitly in terms of local data on $M$ and $G$,
similar to the well-known stratification of $G \backslash M$ by orbit types.

In the case that $G$ is a torus and $M$ is stably almost complex,
an inertia space is implicitly realized as a differentiable stratified space in
\cite{GoHoKnu07}, where the Chen--Ruan orbifold cohomology is extended to this case.
This construction differs from ours in that their inertia space is implicitly defined
as a subquotient of the space $M \times \underline{G}$ where $\underline{G}$ denotes
the group $G$ with the \emph{discrete} topology; the inertia space then appears as the
disjoint union of an infinite family of quotients of
$G$-invariant submanifolds of $M$.  Our construction considers the inertia space as a subquotient of
the manifold $M \times G$ where $G$ is given its usual topology as a Lie group.  Hence, while these
inertia spaces are the same as sets, the topology of our inertia space does not coincide
with that of \cite{GoHoKnu07}.

This paper is organized as follows.  In Section \ref{sec-Prelim}, we review the notions
of differentiable spaces and differentiable stratified spaces.
In Section \ref{sec-StrucSheaf}, we define the inertia
space as well as the structure sheaf with which it is a differentiable space.
In the same section, we also study the local properties of the inertia space and
in particular demonstrate that it is locally
contractible and triangulable.  In Section \ref{sec-StratExamples}, we explicitly describe the
stratification of the inertia space; we give several examples of the main construction before
proving the corresponding main result, Theorem \ref{thrm-StratConnected}.
In Section \ref{sec-deRham}, we prove Theorem \ref{thrm-deRhamInertia},
a de Rham Theorem for the inertia space.

\section*{Acknowledgments}
The first author would like to thank the University of Florence for hospitality during the
completion of this manuscript.
The second author kindly acknowledges NSF support under award DMS 1105670.
The third author would like to thank the University of Colorado at Boulder for its hospitality
and acknowledges support from a Rhodes College Faculty Development Endowment Grant
and a grant to Rhodes College from the Andrew W. Mellon Foundation.


%
%

\section{Preliminaries}
\label{sec-Prelim}
In this section, we recall the definitions of differentiable spaces and
differentiable stratified spaces. Hereby, we use the notion of a differentiable
space as originally introduced by \textsc{Spallek}
\cite{SpallekDR, SpallekGlaettungDR, SpallekDiffForms, SpallekDiffFormsII},
and follow the exposition and notation of \cite{NGonzalezSanchoBook}; see also
\cite{BierstoneLifting, BierstoneOrbitSpace, BredonBook, PflaumBook} for
more details on stratified spaces.

Recall that a \emph{locally $\R$-ringed space} $(X, \mathcal{O})$ consists of a
topological space $X$ equipped with a sheaf $\mathcal{O}$ of $\R$-algebras
such that at each point $x \in X$ the stalk $\mathcal{O}_x$ is
a local ring. Note that by  definition  there then exists
for each point $x \in X$ and open neighborhood $U$ of $x$ an evaluation
map $e_x : \mathcal{O} \rightarrow \R$, $f \mapsto f(x)$.
A \emph{morphism} of locally $\R$-ringed spaces
$(X, \mathcal{O}) \rightarrow (Y,\mathcal{Q})$
is a pair $(f,F)$, where $f:X \rightarrow Y$ is a continuous map
and $F: \mathcal{Q} \rightarrow f_* \mathcal{O}$ a morphism of sheaves
over $Y$ such that for each $x\in X$ the induced map
on the stalks $F_x: \mathcal{Q}_{f(x)} \rightarrow \mathcal{O}_x$
is a local ring homomorphism. Obviously, locally $\R$-ringed spaces
with their morphisms form a category.

A locally $\R$-ringed space is called an
\emph{affine differentiable space}, if for some $n\in \N^*$
there is a closed ideal
$\mathfrak{a} \subseteq \calC^\infty(\R^n)$ such that
$(X, \mathcal{O})$ is isomorphic as a locally $\R$-ringed space
to $(\Specr (A),\mathcal{A})$. Here, $A$ denotes the differentiable algebra
$\calC^\infty(\R^n)/\mathfrak{a}$, $\Specr(A)$ is the \emph{real spectrum} of $A$,
i.e.~the collection of all continuous $\R$-algebra homomorphisms
$A \to \R$ equipped with the Gelfand topology,
and $\mathcal{A}$ is the structure sheaf on $\Specr (A)$, i.e.~the sheaf associated to the
presheaf $U \mapsto A_U$, where $U$ runs through the open sets of $\Specr (A)$ and
$A_U$ is the localization of $A$ over $U$.
A locally ringed space $(X, \mathcal{O})$ is a \emph{differentiable space}, if
for each point $x \in X$ there is an open neighborhood $U$ such that the restriction
$(U, \mathcal{O}_{|U})$ is an affine differentiable space.
If in addition  the map
\[
 \mathcal O (U) \mapsto \calC (U), \: f \mapsto \hat{f} := \big( U \ni x \mapsto f(x) \in \R \big)
\]
is injective  for each open $U\subseteq X$,
one calls  $(X, \mathcal{O})$ a \emph{reduced differentiable space}.
A reduced differentiable space is called a \emph{smooth differentiable space},
if for each point $x \in X$ there is an open neighborhood $U$ such that the restriction
$(U, \mathcal{O}_{|U})$ is isomorphic as a locally $\R$-ringed space to some
$(\R^n, \calC^\infty_{\R^n})$, where $\calC^\infty_{\R^n}$ denotes the sheaf of
smooth functions on $\R^n$.
Examples of reduced differentiable
spaces include smooth manifolds, the orbit space of a proper smooth action of a
Lie group on a manifold or more generally of a proper Lie groupoid \cite{PflPosTanGOSPLG},
algebraic varieties, and symplectically reduced spaces.

Let now $Y \subset X$ be a locally closed subspace of a differentiable
space $(X, \mathcal{O})$. Then, $Y$ carries a natural structure sheaf $\mathcal{O}_{|Y}$
such that $(Y, \mathcal{O}_{|Y})$ becomes a reduced differentiable space. More precisely,
if $V \subseteq Y$ is relatively open and $U\subseteq X$ open with $V= U\cap Y$, the algebra
of all $f\in \calC (V)$ such that there is an  $F \in \mathcal{O} (U)$ with $f(x) =F(x)$ for all
$x\in V $ is independant of the particular choice of $U$ and coincides by definition with
$\mathcal{O}_{|Y} (U\cap Y)$.  Note that in general, the restricted sheaf $\mathcal{O}_{|Y}$
coincides with the pullback sheaf $i^* \mathcal{O}$ for the embedding $i : Y \hookrightarrow X$
only if $Y$ is open in $X$.

Suppose now that $(X, \mathcal{O})$ is a reduced differentiable space which in addition
carries a stratification in the sense of Mather \cite{MatSM}, see also
\cite[Chap.~1]{PflaumBook}.
Then every stratum $S$ of $X$ is locally closed, hence one obtains for every stratum
$S$ the restricted sheaf $\mathcal{O}_{|S}$. Denote by $\calC^\infty_S$
the sheaf of smooth functions on the smooth manifold $S$.
We say that $(X, \mathcal{O})$ is a \emph{differentiable
stratified space}, if for each stratum $S$ of $X$, the sheaves $\mathcal{O}_{|S}$
and $\calC^\infty_S$ coincide.
Note that this notion of a differentiable stratified space is equivalent
to the notion of a stratified space with $\calC^\infty$-structure as defined in
\cite[Sec.~1.3]{PflaumBook}.  In particular, given an affine set
$U$ (i.e. an open subset $U \subseteq X$ such that $(U, \mathcal{O}_{|U})$ is an affine
differentiable space), an isomorphism of $(U, \mathcal{O}_{|U})$ with
$\mbox{Spec}_r(\calC^\infty(\R^n)/\mathfrak{a})$ defines a singular chart
for $U$ in the sense of \cite[Sec.~1.3]{PflaumBook}.
Often, we denote the structure sheaf of a reduced differentiable space $X$ or a
differentiable stratified space $X$ by $\calC^\infty_X$.

To give an example of a differentiable stratified space,
consider a Lie group $G$ acting properly on a smooth manifold $M$.
Let $\varrho : M \to G\backslash M$ be the quotient map.
It is well known (cf.~e.g.~\cite[Sec.~4.3]{PflaumBook} or
\cite[Sec.~2.7]{duistermaatkolk})
that the orbit space $G\backslash M$ is stratified by orbit types.  Specifically,
let $G_x \leq G$ denote the isotropy group of a point $x \in M$, let $(G_x)$
denote the $G$-conjugacy class of $G_x$, and let $M_{(G_x)}$ denote the collection of
$y \in M$ such that $G_y$ is conjugate to $G_x$.  Then the stratification of $G\backslash M$
by orbit types is given by assigning to each $x \in M$ the germ of the set
$G\backslash M_{(G_x)}$.
Moreover, by \cite[Thm.~4.4.6]{PflaumBook}, the orbit space
carries a canonical differentiable structure which is compatible with the stratification by
orbit types. In other words, $G \backslash M$ thus becomes a
differentiable stratified space.  The structure sheaf $\calC^\infty_{G\backslash M}$
is given by assigning to an open subset $U$ of $G\backslash M$ the $\R$-algebra of
continuous functions on $U$ which pull back to smooth $G$-invariant functions on
$\varrho^{-1}(U)$, i.e.
\[
    \calC^\infty_{G\backslash M} (U)
    := \big\{ f \in \calC (U) \mid
    f \circ \varrho_{|\varrho^{-1}(U)} \in \calC^\infty(\varrho^{-1}(U))^G
    \big\} .
\]


%
%

\section{The Inertia Space of a Proper Lie Groupoid}
\label{sec-StrucSheaf}

Recall that by a \emph{groupoid}
one understands a small category $\mathsf{G}$ such that all arrows are
invertible, cf.~\cite{moerdijkmrcun}. Denote by $\sfG_0$ the set
of objects and by $\sfG_1$ the set of arrows of a groupoid $\sfG$.
The source (resp.~target) map will then be denoted by
$s: \sfG_1 \rightarrow \sfG_0$ (resp.~$t: \sfG_1 \rightarrow \sfG_0$),
the unit map by $u:\sfG_0 \rightarrow \sfG_1$, the inversion by
$i:\sfG_1 \rightarrow \sfG_1$, and finally the composition map
by $m :\sfG \times_{\sfG_0} \sfG_1 \rightarrow \sfG_1$. If $\sfG_1$ and $\sfG_0$
are both topological spaces, and all structure maps continuous, the groupoid is
called a \emph{topological groupoid}. If in addition, $\sfG_1$ and $\sfG_0$
are smooth differentiable spaces, all structure maps are smooth maps,
and $s$ and $t$ are both submersions,  $\sfG$ is called a \emph{Lie groupoid}.
Note that the arrow set of a Lie groupoid in general need not be a Hausdorff
topological space.

If $\sfG$ is a topological groupoid, and both $s$ and $t$ are local
homeomorphisms, $\sfG$ is called an \emph{\'etale groupoid}, in case
the pair $(s,t):\sfG_1 \rightarrow |sfG_0 \times \sfG_0$ is a proper map,
one says that $\sfG$ is a \emph{proper groupoid}.

Fundamental examples of  proper Lie groupoids are given by transformation groupoids
$G\ltimes M$, where $G$ is a Lie group which acts properly on a smooth manifold $M$.
The object space of such a transformation groupoid is given by $(G\ltimes M)_0:= M$,
the arrow space by   $(G\ltimes M)_1:= G \times M$, and the structure maps
are defined as follows:
\begin{displaymath}
  \begin{split}
    s & : (G\ltimes M)_1 \rightarrow (G\ltimes M)_0 , \: (g,p) \mapsto p , \\
    t & : (G\ltimes M)_1 \rightarrow (G\ltimes M)_0 , \: (g,p) \mapsto g p , \\
    u & : (G\ltimes M)_0 \rightarrow (G\ltimes M)_1 , \: p \mapsto (e,p) , \\
    i & :(G\ltimes M)_1 \rightarrow (G\ltimes M)_1 , \: (g,p) \mapsto (g^{-1},p),
    \quad \text{and} \\
    m & :  (G\ltimes M)_1 \times_{ (G\ltimes M)_0}  (G\ltimes M)_1\rightarrow (G\ltimes M)_1 , \:
    \big( (g,hp) , (h,p) \big) \mapsto (gh,p) .
  \end{split}
\end{displaymath}

Let now $\sfG$ be an arbitrary proper Lie groupoid. One then defines the \emph{loop space}
of $\sfG$ as the subspace
\begin{equation}
  \label{eq:DefLoopSp}
  \sfB_0:= \big\{ k \in \sfG_1 \mid s(k) = t(k) \big\} \: .
\end{equation}
Sometimes, we denote the loop space also by $\inertia{(\sfG_0)}$.
The groupoid $\sfG$ acts on the loop space in the following way:
\begin{displaymath}
  \sfG_1 \times_{\sfG_0}\sfB_0 \rightarrow \sfB_0, \: (g, k) \mapsto g\, k\, g^{-1} \: .
\end{displaymath}
We can now define:
\begin{definition}[cf.~\cite{LupUri}]
 Let $\sfG$ be a proper Lie groupoid, and $\sfB_0$ its loop space. The
 action groupoid $\sfG\ltimes \sfB_0$ then is called the \emph{inertia groupoid}
 of $\sfG$. It will be denoted by $\inertia{\sfG}$. The quotient space
 $\sfG\backslash \inertia{\sfG}$ will be called the \emph{inertia space} of the
 groupoid $\sfG$. If $X$ denotes the orbit space $\sfG\backslash \sfG_0$,
 we sometimes write (by slight abuse of notation) $\inertia{X} $ for the inertia
 space of $\sfG$.
\end{definition}

\begin{remark}
  The loop space $\sfB_0$ is a closed subset of the smooth manifold $\sfG_1$, hence inherits the structure
  of a differentiable space. Moreover, $\sfB_0$ is locally semialgebraic, hence possesses
  a minimal Whitney B stratification, and a triangulation subordinate to it.
  The inertia  space $\inertia{X}$ inherits these properties from the loop space as well.
  We will elaborate on this in a forthcoming publication.
\end{remark}

\begin{remark}
 The inertia space $\inertia{X}$ depends in fact only on the Morita equivalence class
 of the proper Lie groupoid $\sfG$. So if one thinks of $X$ as a topological space
 together with a Morita equivalence class of Lie groupoids having $X$ as
 orbit space, the notation  $\inertia{X}$ is fully justified.
\end{remark}

Let us now describe the inertia space in the particular situation, where the underlying
proper Lie groupoid is a transformation groupoid $G \ltimes M$ with $G \times M \rightarrow M$
a proper Lie group  action. The loop space $\sfB_0$ then is given as the closed subspace
\[
 \inertia M := \sfB_0 := \big\{ (k,x) \in G \times M \mid kx = x \big\}
\]
of $G\times M$. Moreover, $G$ acts on $G\times M$ by
\[
   G \times (G \times M) \rightarrow (G\times M) , \: (g,(k,x)) \mapsto g (k,x ) := (gkg^{-1} , x ) .
\]
This action leaves $\inertia{M}$ invariant. The inertia space of $G \ltimes M$ now coincides with
the quotient space $\inertia{X} := G \backslash \inertia{M} $, where $X:= G \backslash M$.
Sometimes, we call $\inertia{X}$ the \emph{inertia space of the $G$-manifold $M$}.

\begin{proposition}
  The inertia space  $\inertia{X}$ of a $G$-manifold $M$ carries a natural and
  uniquely determined structure of a differentiable space such that the
  embedding $\iota: \inertia{X} \hookrightarrow G \backslash (G\times M)$
  becomes a smooth map, where $G \backslash (G\times M)$ carries the unique
  differentiable structure such that the canonical projection
  $\varrho   : G\times M \rightarrow G \backslash (G\times M)$ is smooth.
\end{proposition}

\begin{remark}
\label{rem:ProjInertia}
  In the following, we denote the canonical projection
  $M \rightarrow G\backslash M$ of a $G$-manifold $M$ to its orbit space
  by $\varrho^M$. Instead of $\varrho^{G\times M}$ we often
  write $\varrho$, if no confusion can arise. The restriction  of $\varrho$
  to $\inertia{M}$  will be denoted by $\widehat{\varrho}$, that means
  $\widehat{\varrho} : \inertia{M} \rightarrow \inertia{X}$ is the orbit map
  from the loop space to the inertia space.
\end{remark}

\begin{proof}
   Recall that by \cite[Thm.~11.17]{NGonzalezSanchoBook}, the quotient
   $G\backslash (G \times M)$ is a differentiable space, and that the structure sheaf
   on  $G\backslash (G \times M)$ is uniquely determined by the
   requirement that  the quotient map $\varrho$ is smooth.
   Since $\inertia{M}$ is a closed $G$-invariant subspace of
   $G \times M$, it follows from \cite[Lem.~11.15]{NGonzalezSanchoBook}
   that $\inertia{X}$ is a differentiable space. Again, the structure sheaf
   is uniquely determined
   by the requirement that the $\iota: \inertia{X} \hookrightarrow G \backslash
   (G\times M)$ is smooth
   which according to \cite[Lem.~11.15]{NGonzalezSanchoBook} is the case indeed.
\end{proof}

Let us briefly give a more explicit description of the structure sheaf on the inertia space.
Let $U \subset \inertia{X}$ be open. Then $\calC^\infty_{\inertia{X}}(U)$ is the space of all
$f \in \calC (U)$ such that there exists an open $W \subset G \times M$ and a function
$F\in \calC^\infty (W)$ which have the property that $W \cap \inertia{M} = \varrho^{-1} (U)$
and
that $F_{|W \cap \inertia{M}} = f \circ \varrho_{| W \cap \inertia{M}}$. In other words,
\[
    \calC^\infty_{\inertia{X}}(U)
    \cong \big( \calC^\infty (\varrho^{-1} (U))\big)^G.
\]
Hence, the structure sheaf of $\inertia{X}$ is given by
the restriction of the smooth $G$-invariant functions on $G \times M$ to $\inertia{M}$.

The inertia space $\inertia{X}$ of a $G$-manifold $M$ carries even more structure.
In the following considerations we will explain this in more detail.

\subsection{Reduction to Slices in $M$}
\label{subsec-LocalInertiaSlicesM}

Fix a point $x \in M$, and let $Y_x$ be a slice at $x$ for the $G$-action on $M$.
By a \emph{slice at} $x$ we hereby mean a submanifold $Y_x \subset M$ transversely
to the orbit $Gx$ such that the following conditions are satisfied
(cf.~\cite[II. Theorem 4.4]{BredonBook}):
\begin{enumerate}[(SL1)]
\item $Y_x$ is closed in $GY_x$,
\item $GY_x$ is an open neighborhood of $Gx$,
\item $G_x Y_x = Y_x$, and
\item \label{It:SlInv}
      $g Y_x \cap Y_x \neq \emptyset$ implies $g \in G_x$ .
\end{enumerate}
After possibly shrinking $Y_x$, we can even assume that $Y_x$ is a \emph{linear slice},
which means that
\begin{enumerate}[(SL1)]
\setcounter{enumi} 4
\item  \label{It:LinSl}
       there exists a $G_x$-equivariant diffeomorphism
       $Y_x \rightarrow N_x$ of the slice $Y_x$ onto the normal space $N_x := T_xM/T_xGx$.
\end{enumerate}
Note that we implicitly have used here the fact that $G_x$ acts linearly
on the normal space $N_x$.
After choosing a $G$-invariant  riemannian metric on $M$, the image $\exp (B_x)$ of every
sufficiently small open ball $B_x$ around the origin of $N_x$  under the exponential map
is a linear slice at $x$. We will assume from now on that all slices are linear.
By (SL\ref{It:LinSl}) this implies in particular that there is a $G_x$-equivariant
contraction $[0,1] \times Y_x \rightarrow Y_x$ to the point $x$.
Let us also recall at this point the slice theorem  \cite{Koszul}  which tells
that the map
\begin{equation}
  \label{eq:SlThm}
   \theta : G \times_{G_x} Y_x \longrightarrow GY_x , \:
   [h,y] \mapsto hy
\end{equation}
is a $G$-equivariant diffeomorphism between $G \times_{G_x} Y_x$ and the tube
$GY_x \subseteq M$ about $Gx$.

Let us now examine the loop and the inertia space of the $G$-manifold $G \times_{G_x} Y_x$.
By definition, the loop space is given by
\[
  \inertia{(G \times_{G_x} Y_x)} =
  \big\{ (g, [h,y]) \in G \times (G \times_{G_x} Y_x) \mid g[h,y] = [h,y] \big\} .
\]
In addition, the  map
\[
  \id_G \times \theta : \:  G \times (G \times_{G_x} Y_x) \to G \times M
\]
is a $G \times G$-equivariant diffeomorphism onto a $G\times G$-invariant open neighborhood
of $(e,x)$ in $G \times M$. Hence the restriction
\[
  \inertia{\theta} := (\id_G \times\theta )_{|\inertia{(G \times_{G_x} Y_x)} } : \:
  \inertia{(G \times_{G_x} Y_x)} \rightarrow (G \times GY_x ) \cap \inertia{M} \subseteq \inertia{M}
\]
becomes a $G$-equivariant homeomorphism onto the $G$-invariant neighborhood
$(G \times GY_x ) \cap \inertia{M}$ of $(e,x)$ in $\inertia{M}$.
Moreover, it follows that $\inertia{\theta}$ is an isomorphism between the
differentiable spaces $\inertia{(G \times_{G_x} Y_x)}$ and $(G \times GY_x) \cap \inertia{M}$
by \cite[Lem.~11.15]{NGonzalezSanchoBook}.

Next let us consider the loop space
$\inertia{Y_x} := \{ (h,y) \in G_x  \times Y_x \mid hy = y \}$
of the $G_x$-manifold $Y_x$. Then we have the following result,
which provides a local picture of the inertia space of a $G$-manifold $M$.

\begin{proposition}
\label{prop-LocContractRecutionMSlices}
Let $Y_x$ be a slice at the point $x$ of a $G$-manifold $M$.
Then the inertia space $\inertia{(G_x\backslash Y_x)}$ of the $G_x$-manifold $Y_x$ is
isomorphic as a differentiable space to the open neighborhood $\inertia{(G\backslash GY_x)}$
of the point $G(e,x)$ in the inertia space $\inertia{(G\backslash M)}$.
\end{proposition}
\begin{proof}
Consider the map $\phi : \: \inertia{Y_x}  \rightarrow \inertia{(G \times_{G_x} Y_x)} $
defined as the restriction of the smooth map
\[
 G_x \times Y_x \rightarrow G \times (G \times_{G_x} Y_x), \:
 ( h,y) \mapsto   (h, [e, y]).
\]
to $\inertia{Y_x}$. Obviously, by elementary considerations, $\phi$ is continuous
and injective. Moreover, $\phi$  is a morphism of differentiable spaces since it
is the restriction of a smooth map between manifolds.
Since $\phi$
is equivariant with respect to the canonical embedding $G_x \hookrightarrow G$,
it also induces a continuous map between  the quotients
\[
 \Phi: \:  \inertia{(G_x\backslash Y_x)}   \rightarrow
   \inertia{ (G\backslash \left(G \times_{G_x} Y_x\right))} \cong \inertia{(G\backslash GY_x)}, \:
   G_x( h,y ) \mapsto  G(h, [e, y])   .
\]
Let us show that $\Phi$ is bijective. This will prove the claim.

To show that $\Phi$ is injective, suppose that $(h, y)$ and
$(\overline{h}, \overline{y})$ are elements of $\inertia{Y_x}$ such
that $\Phi(G_x(h,y)) = \Phi(G_x(\overline{h}, \overline{y}))$. Then
$G(h, [e, y]) = G( \overline{h} , [e, \overline{y}])$ so that there is
a $g \in G$ such that $(h, [e, y]) = g(\overline{h} , [e, \overline{y}])$.
Therefore, $[e, y] = [g, \overline{y}]$ and $h =
g\overline{h}g^{-1}$, so that there is an $\widetilde{h} \in G_x$ such
that $(\widetilde{h}^{-1}, \widetilde{h}y) = (g, \overline{y})$. But
this implies that $\widetilde{h}^{-1} = g \in G_x$ and
$y = g\overline{y}$, so that $g(\overline{h} ,\overline{y}) = (h, y)$
with $g \in G_x$.  It follows that $\Phi$ is injective.

To show that $\Phi$ is surjective, let $(k, [g, y])$ be an arbitrary
element of the loop space $\inertia{(G \times_{G_x} Y_x)}$ which means that
$k[g, y] = [g, y]$. Then $g^{-1}kg[e,y] = [e, y]$, so that by (SL\ref{It:SlInv})
$g^{-1}kg \in G_x$.
Since $g^{-1}(k, [g, y]) = ( g^{-1}kg, [e, y])$, it follows that
$\Phi(G_x(g^{-1}kg , y )) = G(k, [g, y])$, and $\Phi$ is surjective.

Finally, we claim that $\Phi$  is even an isomorphism between differentiable spaces.
To this end note first that for all
$f \in \calC^\infty \big( \inertia{ (G\backslash \left(G \times_{G_x} Y_x\right))}\big)$
the pullback $\Phi^* (f)$ is a smooth function on  $\inertia{G\backslash Y_x}$,
since
\[
 \Phi^* (f) \circ \varrho^{Y_x}_{|\inertia{Y_x}} =
 f \circ \varrho^{G \times_{G_x} Y_x}_{|\inertia{(G \times_{G_x} Y_x)}}  \circ \phi,
\]
where we have used the notation as explained in Remark \ref{rem:ProjInertia}.
By surjectivity of $\Phi$, the pullback
$\Phi^* : \: \calC^\infty \big( \inertia{ (G\backslash \left( G \times_{G_x} Y_x\right))}
\big) \rightarrow \calC^\infty \big( \inertia{(G_x\backslash Y_x)} \big)$ is injective.

To show that $\Phi^*$ is surjective, let
$h \in \calC^\infty \big( \inertia{ (G_x \backslash Y_x)} \big)$.
Since $\inertia{Y_x}$ is a closed differentiable subspace of $G \times Y_x $, there exists
a smooth $H$ on $G \times Y_x $, such that
$H_{| \inertia{Y_x}} = h \circ \varrho^{Y_x}_{|\inertia{Y_x}}$.
By possibly averaging over $G_x$ one can even assume that  $H$ is $G_x$-invariant.
With this, we define
$\tilde f\co G \times (G \times_{G_x} Y_x) \to \R$ by setting
$\tilde f (k, [g,y]) = H (g^{-1}kg,y )$. By $G_x$-invariance, $\tilde f$ is well-defined
and smooth. Moreover, $\tilde f$ is $G$-invariant, hence on has
\[
  \tilde f_{|\inertia{G \times_{G_x} Y_x}} = f \circ \varrho^{G \times_{G_x} Y_x}_{|\inertia{(G \times_{G_x} Y_x)}}
\]
for some smooth $f: \inertia{(G \times_{G_x} Y_x)} \rightarrow \R$.
By construction it is clear that then $\Phi^* (f) = h$, hence $\Phi^*$ is surjective.
This proves that $\Phi$ is even an isomorphism of differentiable spaces.
\end{proof}

In the preceding proposition, one can choose for each $x\in M$ the slice $Y_x$ to be equivariantly
isomorphic to a ball $B_x$ in some finite dimensional orthogonal $G_x$-representation space $V_x$.
Moreover, since every compact Lie group has a linear faithful representation, we can assume
that $G_x$ is (represented as) a compact real algebraic group. The loop space
$\inertia{B_x} = \big\{ (k,x) \in G_x \times B_x  \mid kx =x \big\}$ then is a  semi-algebraic set in
$G_x \times V_x$. Since by the Proposition the inertia space
$\inertia{(G\backslash M)}$ has a
locally finite cover by open subsets such that each of the elements of the cover is isomorphic as a
differentiable space  to some inertia space $\inertia{(G_x \backslash B_x)}$ the inertia space of the $G$-manifold $M$
is locally semi-algebraic in the sense of Delfs--Knebusch \cite[Chap.~I]{DelKneLSS}.

\begin{corollary}
  The inertia space $\inertia{(G\backslash M)}$ of a $G$-manifold $M$ is locally semi-algebraic.
\end{corollary}

A semi-algebraic set $X$ has a minimal $\calC^\infty$-Whitney stratification according to
\cite[Thm.~4.9 \& p.~210]{MatSM}. By definition of a stratification as in \cite[Sec.~I.2]{MatSM}
it follows that every locally semi-algebraic space carries a minimal $\calC^\infty$-Whitney stratification
compatible with the differentiable structure.
By the Corollary, the inertia space of a $G$-manifold $M$ therefore possesses a minimal $\calC^\infty$-Whitney
stratification as well, and becomes a differentiable stratified space.
We will call this stratification the \emph{canonical stratification} or the
\emph{minimal Whitney stratification} of the inertia space.
Since the canonical stratification satisfies Whitney's condition B, there even exists a system of smooth
control data for the canonical stratification (cf.~\cite{MatNTS} and \cite[Thm.~3.6.9]{PflaumBook}). According
to \cite[Sec.~5]{gor} \cite[Cor.~3.7]{verona} there  exists a triangulation of the inertia space subordinate
to the canonical stratification.
Hence, the inertia space is triangulable and in particular locally contractible.
We thus obtain the following result.

\begin{theorem}
\label{thm:inertiatopprop}
  The inertia space $\inertia{X}$ of a $G$-manifold $M$ is in a canonical way a locally compact locally
  contractible differentiable stratified space.
  Its canonical stratification satisfies Whitney's condition B and is minimal among
  the stratifications of $\inertia{X}$ with this property. Moreover, there exists a triangulation
  of the inertia space subordinate to the canonical stratification.
\end{theorem}

That the inertia space is locally contractible can be shown directly as well. The main point
hereby lies in the following result which also will be needed later to prove a de Rham Theorem
for inertia spaces.

\begin{proposition}
\label{Prop-LocalInertiaScalarMult}
  Let $M$ and $G$ be as above, and consider a point $(h,x)$ in the loop space $\inertia{M}$.
  Then there is a linear slice $V_{(h,x)}$ at $(h,x)$ such that the action of
  scalars $t \in [0, 1]$ on $V_{(h,x)}$ leaves  the set $V_{(h,x)} \cap \inertia{M}$
  invariant. The linearization  $V_{(h,x)} \hookrightarrow B_{(h,x)}\subseteq N_{(h,x)}$ from
  $V_{(h,x)}$ to an open  convex neighborhood of the origin in the normal space
  $N_{(h,x)} := T_{(h,x)} (G \times M) / T_{(h,x)} \big(G (h,x)\big) $ can hereby be chosen
  as the inverse of the restriction $\exp_{|B_{(h,x)}}$ of the exponential map
  corresponding to an appropriate $G$-invariant riemannian metric on $G\times M$.
\end{proposition}
\begin{proof}
Choose a $G$-invariant riemannian metric on $M$, a bi-invariant riemannian metric on $G$
(cf.~\cite[Prop.~2.5.2 and Sec.~3.1]{duistermaatkolk}),
and let $G \times M $ carry the product metric. Recall that under these assumptions,
the exponential map $T_{(h,x)}(G \times M) \to G \times M $ decomposes into a product
$\exp_h^G \times \exp_x^M$. Moreover, the
riemannian exponential map $\exp^G$ on $TG$ then coincides with  the map
\[
  TG  \cong G \times \frak{g} \rightarrow G , \: (g,\xi) \mapsto g \, e^\xi .
\]
Hereby, $e^\xi$ denotes the exponential map on the Lie algebra, and the isomorphism between
$G \times \frak{g}$ and $TG$ is given by
$(g,\xi)\mapsto {(L_g)}_\ast  \xi$, where $L_g : G \rightarrow G$ denotes the left action by $g$.
Now let $B_{(h,x)}$ denote a sufficiently small open ball around the origin of the normal space
$N_{(h,x)}$ so that the exponential map is injective on $B_{(h,x)}$, and let $V_{(h,x)} := \exp\big( B_{(h,x)} \big)$.

Note that every element $g$ of the isotropy group $H := G_{(h,x)}$
commutes with $H$. Hence, for $(k,y) \in V_{(h,x)}$ one has
$h = ghg^{-1} \in gH_{(k,y)}g^{-1}$ if and only if $h \in H_{(k,y)}$,
and the same holds for each connected component of $H_{(k,y)}$.
Therefore, if $(H_1), \ldots, (H_s)$ denotes the collection of isotropy types for the linear
$H$-action on $V_{(h,x)}$, which is finite by \cite[Lem.~4.3.6]{PflaumBook},
we can assume they are ordered in such a way that $h \in gH_i g^{-1}$ for each $g \in H$ and
$i = 1, \ldots, r$, and $h \not\in gH_i g^{-1}$ for each $g \in H$ and
$i = r+1, \ldots, s$.  Moreover, for $i = 1, \ldots, r$, let
$H_i^h$ denote the connected component of $H_i$ containing $h$, which implies that
$h \in gH_i^hg^{-1}$ for each $g \in H$.

With this, we define
\[
    C =
    \left(\bigcup\limits_{i=1}^r \;\; \bigcup\limits_{g \in H} gH_ig^{-1}
    \smallsetminus gH_i^hg^{-1}\right)
    \cup
    \left(\bigcup\limits_{i=r+1}^s \;\; \bigcup\limits_{g \in H} gH_ig^{-1}\right).
\]
That is, $C$ is the union of all conjugates of isotropy groups not containing $h$ as well
as, for each isotropy group containing $h$, all conjugates of the connected components not
containing $h$.
Since  the quotient map $H \to Ad_H\backslash H$ is closed by \cite[Prop.~3.6]{tomdieck}
it follows immediately that  $C$ is closed in $H$. By construction,  $C$ is also $H$-invariant.
This implies that $G \smallsetminus C$ is an $Ad_H$-invariant open neighborhood of $h$ in $H$.
Hence there exists a connected open and $Ad_H$-invariant neighborhood $O_h$ of $h$ in
$G \smallsetminus C$ small enough to be contained in a logarithmic chart around $h$ in the Lie
algebra $\mathfrak{h}$ of $H$, see \cite[Thm.~1.6.3]{duistermaatkolk}.
Therefore, $O_h \times M$ is an $H$-invariant open neighborhood
of $(h,x)$ in $ G \times M$.  With this, we may shrink $V_{(h,x)}$ to assume that
$V_{(h,x)} \subseteq O_h \times M $.

Now, suppose $(k,y) \in V_{(h,x)} \cap \inertia{M}$ so that $ky =y$ and then clearly
$k(k,y) = (k,y)$.  By property (SL\ref{It:SlInv}) of the slice $V_{(h,x)}$, it follows that
$k \in H$. Since $k \in O_h$ we have that $k$ is not contained in any $H$-conjugate of $H_i$ for
$i = r+1, \ldots , s$ so that $H_{(k,y)}$ is conjugate to $H_i$ for
some $i \leq r$.  Moreover, $O_h$ is connected and does not
intersect $C$ so that $k$ is contained in the same connected
component of $H_{(k,y)}$ as $h$.  Hence, by using logarithmic
coordinates near $h$, we may express $k = he^\xi$ for some $\xi \in
\mathfrak{h}_{(k,y)}$, the Lie algebra of $H_{(k,y)}$.
Additionally, $he^{t\xi} \in H_{(k,y)}$ for $t \in [0,1]$, so that
$he^{t\xi}(k,y) = (k,y)$. Next, let $w \in T_xM$ such that $\exp^M_x (w) =y$.
Then we have $\exp_{(h,x)}(\xi,w) = (k,y)$, and $(\xi,w) \in N_{(h,x)}$. Moreover, we get
\[
   \exp_{(h,x)} \big(t(\xi,w) \big) =  \exp_{(h,x)} \big( (t\xi, tw) \big) = \big( he^{t\xi}  , y(t) \big),
\]
where $y(t) = \exp^M_x (tw)$. However, since the action of $H$ on the normal space
$N_{(h,x)}$ is linear, and since $he^{t\xi}(k,y) = (k,y)$, it follows
that $he^{t\xi}\big( he^{t\xi}, y(t)\big) = \big( he^{t\xi}, y(t)\big)$ for $t \in [0,1]$.
This of course implies that $he^{t\xi}y(t) = y(t)$ so that $\big( he^{t\xi}, y(t)\big) \in \inertia{M}$,
proving the claim.
\end{proof}

\begin{corollary}
\label{cor-LocalInertiaContractible} The inertia space
$\inertia{(G\backslash M)}$ of a compact Lie group action is locally contractible.
\end{corollary}
\begin{proof}
Let $(h,x) \in \inertia{M}$, $H = G_{(x,h)}$, and choose a linear
slice $V_{(h,x)}$ as in the preceding Proposition. Then $GV_{(h,x)}$ is an open
$G$-invariant neighborhood of $(h,x)$ in $G\times M $. Define the map
$\mathcal{H} \co G \times_H V_{(h,x)} \times [0,1] \to G \times_H V_{(h,x)}$ by
$\mathcal{H}([g, (k,y)], t) = [g, (1-t)(k,y)]$.
Then $\mathcal{H}$ is a $G$-equivariant deformation retraction of
$G \times_H V_{(h,x)}$ onto $G \times_H \{ (h,x) \}$ which induces a
retraction in the quotient onto the single orbit $G(h,x)$. Moreover,
by the preceding Proposition, the map $\mathcal{H}$
restricts to a $G$-invariant retraction of $(G \times_H V_{(h,x)}) \cap \inertia{M}$
onto a single orbit.
\end{proof}



%
%

\section{The orbit Cartan type Stratification}
\label{sec-StratExamples}

In this section, we present the explicit stratification of the
inertia space $\inertia{X}$.  We give the definition of this stratification
in Subsection \ref{subsec-ExamplesStratDescription} and state our first main result,
Theorem \ref{thrm-StratConnected}.  Before turning to the proof of
Theorem \ref{thrm-StratConnected} in Subsection \ref{subsec-StratConnectedProof},
we first give several examples of the stratifications in Subsection \ref{subsec-Examples} and
establish some useful results for actions of abelian groups
in Subsection \ref{subsec-StratConnectedTorus}.

Let $G^\circ$ denote the connected component of the identity of $G$.
Recall that a \emph{Cartan subgroup} $\cartan$ of $G$ is a closed topologically cyclic subgroup that has finite
index in its normalizer $N_G(\cartan)$ (cf.~\cite[IV.~Def.~4.1]{BroeckertomDieck}, see also \cite{Segal}).
If $g \in G$, then by
\cite[IV.~ Prop.~4.2]{BroeckertomDieck}, there is a Cartan
subgroup $\cartan$ of $G$ such that $\cartan/\cartan^\circ$ is generated
by $g\cartan^\circ$.
We will say that such a $\cartan$ is a \emph{Cartan subgroup associated to $g$}.
If $g \in G^\circ$, then $\cartan$ is a maximal torus of $G^\circ$ containing $g$;
in general, $\cartan$ is isomorphic to the product of a torus and a finite
cyclic group.  We will make frequent use of \cite[IV.~Prop.~4.6]{BroeckertomDieck},
which states that the homomorphism
\begin{equation}
\label{eq-StratConnectedCartanConjugate}
\begin{array}{rcl}
    \cartan/\cartan^\circ
        &\longrightarrow&           G/G^\circ           \\
    t\cartan^\circ
        &\longmapsto&               tG^\circ
\end{array}
\end{equation}
defines a correspondence between Cartan subgroups of $G$ and cyclic subgroups of
$G/G^\circ$ that induces a bijection on conjugacy classes.  That is, given $g, h \in G$
and Cartan subgroups $\cartan_g$ and $\cartan_h$ associated to $g$ and $h$, respectively,
$\cartan_g$ and $\cartan_h$ are conjugate in $G$ if and only if
$\langle gG^\circ \rangle \leq G/G^\circ$ and $\langle hG^\circ \rangle \leq G/G^\circ$
are conjugate in $G/G^\circ$.  For $g, h \in G^\circ$, this corresponds to the well-known
fact that all maximal tori in a compact connected Lie group are conjugate; see e.g.
\cite[Thm.~3.7.1]{duistermaatkolk}.

\subsection{Definition of the Stratification}
\label{subsec-ExamplesStratDescription}

Let $(h, x) \in \inertia{M}$ and let $H$ denote the isotropy
group of $(h, x)$ with respect to the $G$-action on $G \times M$.
Then $H = G_x \cap Z_G(h) = Z_{G_x}(h)$ where $Z_G(h)$ denotes the
centralizer of $h$ in $G$ and $G_x$ denotes the isotropy group of
$x$ with respect to the $G$-action on $M$.  Let $\cartan_{(h,x)}$
be a Cartan subgroup of $H$ associated to $h$; note that if
$G_x$ is connected, we have by \cite[Thm.~3.3.1 (i)]{duistermaatkolk}
that $h \in (Z_{G_x}(h))^\circ = H^\circ$, so that $\cartan_{(h, x)}$
is a maximal torus of $H^\circ$ containing $h$.
Choose a slice $V_{(h,x)}$ at $(h, x)$ for the $G$-action on $G \times M$, and
define an equivalence relation $\sim$ on $\cartan_{(h,x)}$
by $s \sim t$ if $(GV_{(h,x)})^s = (GV_{(h,x)})^t$.
Let $\cartan_{(h,x)}^\ast$ denote the connected component of the
$\sim$ class $[h]_{}$ containing $h$.

We define a stratification of $\inertia{M}$ by assigning to the point $(h,x) \in \inertia{M}$
the germ
\begin{equation}
\label{eq-StratTildeMDef}
    \mathcal{S}_{(h,x)}
    =
    \left[G\left( V_{(h,x)}^H \cap
    \left(\cartan_{(h,x)}^\ast \times M\right)\right)\right]_{(h,x)}.
\end{equation}
After applying the quotient map $\widehat{\varrho}: \inertia{M} \to \inertia{X}$,
we can similarly define a stratification of $\inertia{X}$
by assigning to the orbit $G(h,x)$ the germ
\begin{equation}
\label{eq-StratInertDef}
    \mathcal{R}_{G(h,x)}
    =
    \left[\widehat{\varrho}\left(G\left( V_{(h,x)}^H \cap
    \left(\cartan_{(h,x)}^\ast\times G\right)\right)\right)\right]_{G(h,x)}.
\end{equation}

The following result shows that $\mathcal{S}$ and $\mathcal{R}$ are stratifications of the loop space and
inertia space, indeed. We call these the stratifications by \emph{orbit Cartan type}.

\begin{theorem}
\label{thrm-StratConnected}
Let $G$ be a compact Lie group and let $M$ be a smooth $G$-manifold.
Then Equation \eqref{eq-StratTildeMDef} defines a Whitney stratification of
$\inertia{M}$ with respect to which $\inertia{M}$ is a differentiable stratified space.
Moreover, this stratification induces a stratification on $\inertia{X}$
by Equation \eqref{eq-StratInertDef} with respect to which
$\inertia{X}$ is a differentiable stratified space fulfilling Whitney's condition B.
\end{theorem}

An immediate consequence of this result and Theorem \ref{thm:inertiatopprop} is the following.

\begin{corollary}
  The orbit Cartan type stratification is in general finer than the canonical stratification
  of the inertia space $\inertia{X}$.
  Moreover, there exists a triangulation
  of the inertia space subordinate to the orbit Cartan type stratification.
\end{corollary}

\begin{remark}
\label{rem-StratGlobalDef}
Assume in addition that $M$ itself is partitioned into a finite number of $G$-
and $\cartan$-isotropy types for any Cartan subgroup $\cartan$ of $G$.
The definition of $\cartan_{(h,x)}^\ast$ above can be modified by saying that
$s \sim t$ if $(G \times M)^s = (G \times M)^t$. The modified definitions of
$\mathcal{S}_{(h,x)}$ and $\mathcal{R}_{G(h,x)}$ also result in Whitney stratifications
of $\inertia{M}$ and $\inertia{X}$, respectively.  The proof of this fact
is identical to the proof of Theorem \ref{thrm-StratConnected} below with minor simplifications.
The modified stratifications are generally finer and depend on global data in $G \times M$,
though they can be easier to compute in examples.
\end{remark}

Before we prove Theorem \ref{thrm-StratConnected}, we provide several examples which illustrate our definition.

\newpage

\subsection{Examples of the Stratification}
\label{subsec-Examples}

\subsubsection{Cases Where $\inertia{M}$ is Smooth}
\label{subsubsec-ExamplesFree}

Suppose $G$ acts freely on $M$. Then
\[
    \inertia{M} = \{ e \} \times M \subseteq G \times M
\]
is diffeomorphic to $M$.  Each point $(e, x)$ has trivial isotropy group,
and it is easy to see that the stratifications of $\inertia{M}$ and $\inertia{X}$
given by Equations \eqref{eq-StratTildeMDef} and \eqref{eq-StratInertDef} are trivial.
The result in both cases is a smooth manifold with a single stratum, and hence trivially a stratified
differentiable space.

Similarly, suppose $L$ is a (necessarily normal) subgroup of $G$ that acts
trivially on $M$, and suppose $G/L$ acts freely on $M$.  Then
\[
    \inertia{M} = L \times M \subseteq G \times M
\]
is a smooth manifold.  The isotropy types of elements of $\inertia{M}$ correspond to
the isotropy types of $L$ with respect to the $G$-action by conjugation; that is,
elements $(h,x)$ and $(k, y)$ of $\inertia{M}$ have the same isotropy type if and
only if the centralizers $Z_G(h)$ and $Z_G(k)$ are conjugate.  We claim that in this case,
the stratifications of $\inertia{M}$ and $\inertia{X}$ given by
Equations \eqref{eq-StratTildeMDef} and \eqref{eq-StratInertDef} coincide with the
stratifications by $G$-isotropy types.

Choose a slice $V_{(h,x)}$ at $(h,x) \in \inertia{M}$ for the $G$-action on $G \times M$.
By construction, it is clear that $\mathcal{S}_{(h,x)}$ is a subgerm
of the germ of the isotropy type of $(h,x)$ at $(h,x)$.
Let $(\tilde{k},\tilde{y}) \in G V_{(h,x)} \cap ( L \times M )$ be a point in the orbit of this slice
with the same $G$-isotropy type as $(h,x)$.  Then there is a
$\tilde{g} \in G $ such that $(k,y) := \tilde{g}(\tilde{k},\tilde{y}) \in V_{(h,x)}$, and hence
$G_{(k,y)} \leq H = G_{(h,x)}$.  However, as $G_{(k,y)} = \tilde{g}G_{(\tilde{k},\tilde{y})}\tilde{g}^{-1}$
is conjugate to $H$, we have by \cite[Lem.~4.2.9]{PflaumBook}
that $G_{(k,y)} = H$.  Therefore,
$(k,y) \in V_{(h,x)}^H$, which is connected, so that $k$ is in the same connected component of
$H$ as $h$.  It follows that $kH^\circ$ and $hH^\circ$ generate the same subgroup of
$H/H^\circ$, so that by \cite[IV.~Prop.~4.6]{BroeckertomDieck}, Cartan subgroups
$\cartan_{(h,x)}$ and $\cartan_{(k,y)}$ of $H$ associated to $h$ and $k$, respectively, are
conjugate in $H$.  Hence there is a $g \in H$
such that $gkg^{-1} \in \cartan_{(h,x)}$; however, as $g \in H = G_{(k,y)}$,
it follows that $k = gkg^{-1} \in \cartan_{(h,x)}$.
Moreover, as $Z_L(k) = G_{(k,y)} = G_{(h,x)} = Z_L(h)$,
we have that $(G \times M)^k = (G \times M)^h$ so that $(GV_{(h,x)})^h = (GV_{(h,x)})^k$ and
$k \in \cartan_{(h,x)}^\ast$.  We conclude that
$(\tilde{k},\tilde{y}) \in G\left( V_{(h,x)}^H \cap \left(\cartan_{(h,x)}^\ast\times M\right)\right)$,
and hence that $\mathcal{S}_{(h,x)}$ is the germ of the $G$-isotropy type of $(h,x)$ in $M \times G$.

More generally, we have the following.  The proof is an elementary argument applied to slices
for the $G$-action on $M$ using a local section of the fiber bundle $G\to G/J$.

\begin{proposition}
\label{prop-StratConnectedDepthZero}
Suppose the stratification of $M$ by $G$-orbit types has depth zero which means that there is
a $K \leq G$ such that every point has orbit type $(K)$.
Then $\inertia{M}$ is a smooth submanifold of $G \times M$ that is locally diffeomorphic
to $K \times M$.
\end{proposition}

\begin{proof}
To show that $\inertia{M}$ is a differentiable manifold, let $(h, x) \in \inertia{M}$,
and let $Y_x$ be a slice at $x$ for the $G$-action on $M$.  Without loss of generality,
we can assume that $G_x = K$.  Then for each $y \in Y_x$, as $G_y \leq K$
and $G_y$ is conjugate to $K$, it must be 
that $G_y = K$.  Therefore,
$Y_x^{K} = Y_x$, and a neighborhood of $G_x$ in $M$ is diffeomorphic to
\[
    G \times_{K} Y_x  =   G/K \times Y_x
\]
via the map
\begin{displaymath}
    \tau   :\:    G/K \times Y_x  \longrightarrow   M , \:  (gK, y)  \longmapsto   gy.
\end{displaymath}

To prove that $\inertia{M}$ is a differentiable submanifold of $G \times M$,
choose a neighborhood $U$ of $eK$ in $G/K$ small enough so that the fiber bundle $G \to G/K$ admits a
differentiable section on $U$.  Let $\sigma : U \to G$ for $G \to G/K$ be a choice of such a section,
and consider the map
\[
\begin{split}
    \widetilde{\tau} :\; & G \times U \times Y_x   \longrightarrow
    G \times \tau(U \times Y_x) \subseteq G \times M , \\
    & (\widehat{g}, gK, y) \longmapsto
    \big( \sigma(gK) \, \widehat{g} \, \sigma(gK)^{-1}, \sigma(gK)y \big) \: .
\end{split}
\]
Since $U$ is an open neighborhood of $eK$ in $G/K$, we have that $U \times Y_x$
is an open neighborhood of $(eK, x)$ in $G/K \times Y_x$.  Therefore,
$\tau(U \times Y_x)$ is an open neighborhood of $x$ in $M$.
Simple computations demonstrate that $\widetilde{\tau}$ is a diffeomorphism from the neighborhood
$G \times U \times Y_x$ of $(h, eK, x)$ in $G \times G/K \times Y_x$
onto the neighborhood $G\times \tau(U \times Y_x)$ of $(h, x)$ in $G \times M$.  Moreover,
\[
    \widetilde{\tau}(K \times U \times Y_x)
    = (G \times \tau(U \times Y_x)) \cap \inertia{M},
\]
so that $\widetilde{\tau}$ restricts to a diffeomorphism between a neighborhood
of $(h, x)$ in $K \times M$ to a neighborhood of $(h, x)$ in $\inertia{M}$.
\end{proof}

In this case, however, it may happen that
the stratifications of $\inertia{M}$ and $\inertia{X}$
given by Equations \eqref{eq-StratTildeMDef} and \eqref{eq-StratInertDef} are strictly
finer than the respective stratifications by isotropy types.
This is the case, for instance, when $\inertia{X}$ is the inertia space of
$\R^3 \smallsetminus \{ 0 \}$ with its usual
$\SO{3}$-action; see \ref{subsubsec-ExamplesSO3} below.

\subsubsection{Locally Free Actions}
\label{subsubsec-ExamplesLocallyfree}

If the action of $G$ on $M$ is locally free, i.e.~the isotropy group of each $x\in M$ is finite,
then the quotient $X = G\backslash M$ is an orbifold. The corresponding inertia space $\inertia{X}$
then is an orbifold as well and is called the \emph{inertia orbifold} of $X$, see e.g.~\cite{ademleidaruan} or
\cite{PflaumAlgIndex}. Let us briefly sketch this within our framework and let us show that the above defined
stratification of the inertia space $\inertia{X}$  coincides with the orbit
type stratification, if the action of $G$ is locally free.

To this end consider first the case, where $G$ is a finite group.
The loop space $\inertia{M}$ then is the disjoint union $\bigsqcup_{h \in G} \{ h \} \times M^h$ of
smooth manifolds of possibly different dimensions. Choose a $G$-invariant riemannian metric on $M$.
Since $G$ is finite, the linear slice $V_{(h,x)}$ at some point $(h,x) \in \inertia{M}$
can be chosen to be of the form $\{ h \} \times V_x$, where $V_x \subset M^h$ is an open ball around $x$
in $M^h$; note that $M^h$ is totally geodesic in $M$. Denote by $H$ the isotropy group of $(h,x)$, i.e.~let
$H:=Z_{G_x} (h)$. Because under the assumptions made the Cartan subgroups are discrete,
the set germ $\mathcal{S}_{(h,x)}$ at $(h,x) \in \inertia{M}$ from Eq.~\eqref{eq-StratTildeMDef}
coincides with
\[
  [G(V^H_{(h,x)} \cap (\{ h \} \times M))]_{(h,x)} =
   [G(\{ h \} \times V_x^{H})]_{(h,x)} = [\{ h \} \times V_x^{H}]_{(h,x)} \, .
\]
The second equality hereby follows from the fact that for every $g \in Z_G(h)$ and $y \in V_x^H$  with
$gy  \in V_x$  one has $gy \in  V_x^H$, since $g \in H$ by (SL\ref{It:SlInv}).
Observe that the orbit map $\varrho : G \times M \rightarrow  G \backslash (G \times M )$ is injective
on $\{ h \} \times V_x^{H}$ by the slice theorem, hence the set germ  $\mathcal{R}_{G(h,x)}$ at $G(h,x) \in \inertia{M}$
is given by
\[
   \mathcal{R}_{G(h,x)} = [\varrho (\{ h \} \times V_x^{H})]_{G(h,x)} \, .
\]
In other words this means that the stratification by orbit Cartan type of the inertia space $\inertia{X}$
of a finite group action on $M$ is given by the orbit type stratification.

Let us now consider the case where $G$ is a compact Lie group acting locally freely on $M$.
According to Theorem \ref{thm:inertiatopprop}, the inertia space $\inertia{X}$ is a differentiable stratified space
with stratification given by the canonical stratification. Recall that  the canonical stratification is
minimal among all Whitney stratifications of $\inertia{X}$.
Now observe that by Proposition \ref{prop-LocContractRecutionMSlices} the neighborhood
$\inertia{(G \backslash GY_x)}$ of $(e,x)$ in $\inertia{X}$
is isomorphic as a differentiable space to the inertia space $\inertia{(G_x\backslash Y_x)}$,
where $Y_x $ is a slice of $M$ at $x$.  Since $G_x$ is finite, it follows by the above considerations that
the stratification $\mathcal{R}$ of $\inertia{(G_x\backslash Y_x)}$ coincides with the stratification by orbit
types. But the latter is known to be the minimal Whitney stratification, hence $\inertia{(G_x\backslash Y_x)}$
with the orbit type stratification is even isomorphic as a differentiable stratification to
$\inertia{(G \backslash GY_x)}$ with the canonical stratification. Since  $\inertia{X}$ is covered by
the open sets $\inertia{(G \backslash GY_x)}$, $x \in X$, it follows that both the canonical stratification
of $\inertia{X}$ and the stratification by orbit Cartan type coincide with the stratification by orbit type, and
that $\inertia{X}$ is an orbifold, indeed.

\subsubsection{Semifree Actions}
\label{subsubsec-ExamplesSemifree}

Suppose $G$ acts semifreely on $M$ so that there is a
collection $N = M^G$ of submanifolds of $M$ fixed by $G$, and $G$
acts freely on $M \smallsetminus N$.  Then
\[
    \inertia{M}
    =
    [\{ e \}\times (M \smallsetminus N)] \cup (G \times N).
\]
The isotropy group of $(h,x)$ is trivial if $x \notin N$ and is equal to the centralizer
$Z_G(h)$ if $x \in N$.  With respect to the adjoint action $Ad_G$ of $G$ on itself, let
$(K_1), \ldots, (K_m)$ denote the isotropy types of elements of $G$
so that the centralizer of every element of $G$ is conjugate to some
$K_j$.  We assume that $K_m = G$ is the isotropy group of the center of $G$.
For $j = 1, \ldots, m - 1$, we let $Ad_G^j$ denote the set of elements of $G$
with centralizer exactly $K_j$ and let $Ad_G^{(j)}$ denote the set of elements of $G$ with
centralizer conjugate to $K_j$.  For $j = m$, we let $Ad_G^m = Ad_G^{(m)}$
denote the set of \emph{nontrivial} central elements of $G$, which
may be empty.  The sets $Ad_G^j$ and $Ad_G^{(j)}$ are
disjoint unions of smooth submanifolds of $G$ by \cite[Cor.~4.2.8]{PflaumBook}.
Moreover, we have for $h \in G$ that $(G \times M)^h = Z_G(h) \times N$, so that
if $(h,x)$ and $(k, y)$ have the same isotropy group, then $h$ and $k$ have the same
fixed point sets in neighborhoods of the orbits $G(h,x)$ and $G(k,y)$.

Let
\[
    \mathcal{S}_0 :=
    \{ e \} \times (M \smallsetminus N),
\]
and for each $j \in \{1, \ldots, m\}$, let
\[
    \mathcal{S}_j :=
    Ad_G^{(j)} \times N
\]
(which is empty for $j = m$ if $G$ has trivial center), and
\[
    \mathcal{S}_{m + 1} :=
    \{ e \} \times N.
\]
Projection under the quotient map
$\widehat{\varrho} : \inertia{M} \to \inertia{X}$
provides manifolds
\[
    \mathcal{R}_j :=
    \widehat{\varrho}(\mathcal{S}_j).
\]
The decompositions of $\inertia{M}$ and $\inertia{X}$ given by the
connected components of the $\mathcal{S}_j$ and $\mathcal{R}_j$, respectively,
coincide with the stratifications defined in Equations
\eqref{eq-StratTildeMDef} and \eqref{eq-StratInertDef}.

In particular, note that this stratification is strictly finer than the stratification
by orbit types in the case where $G$ has nontrivial center.  In fact, the piece
$Z(G) \times N$, which consists of points of the same isotropy type, must be split
into $\{ e \} \times N$ and $(Z(G)\smallsetminus \{ e \}) \times N$ in order
for the pieces to satisfy the condition of frontier.  The reason for
this is the occurrence of $\{ e \}$ as the isotropy group for points
of the form $(e, x)$ with $x \in M \smallsetminus N$.  Indeed, the closure of
the strata $\mathcal{S}_0 =  \{ e \} \times (M \smallsetminus N)$
is $\{e \}\times M$, and hence cannot contain the entire isotropy type
of points $(e, n)$ with $n \in N$.

As a simple, concrete example, consider the action of the circle $\SO{2}$ on the sphere $\Sp^2$
by rotations about the $z$-axis; this action is semifree with $N = (\Sp^2)^{\SO{2}}$ given by
the north and south poles.  It is easy to see that the isotropy types
\[
    A = \{ (e, x) \mid x \in \Sp^2 \smallsetminus N \}
\]
and
\[
    B = \{ (t, x) \mid t \in S^1, x \in N \}
\]
do not yield a decomposition of $\inertia{\Sp^2}$,
as $\overline{A} \cap B = \{ e \} \times N$.


\subsubsection{Actions of Abelian Groups}
\label{subsubsec-ExamplesAbelian}

Suppose $G$ is abelian, and let $\{ H_i \mid i \in I \}$ be the
(possibly infinite) collection of isotropy groups for the
$G$-action on $M$.  Note that the isotropy group of $(h,x) \in
\inertia{M}$ is equal to $G_x$. For each $x \in M$, let $I_x
\subseteq I$ be the finite subset consisting of all $i$ such that
every neighborhood of $x$ contains points with isotropy group $H_i$.

Choose $(h,x) \in \inertia{M}$, and note that the Cartan subgroup $\cartan_{(h,x)}$
is in this case unique.  For $k \in \cartan_{(h,x)}$, we have $k \sim h$ if and only
if $h$ and $k$ fix the same points in a neighborhood of $x$ in $M$, or equivalently if
and only if $h$ and $k$ are in exactly the same isotropy groups $H_i$ for $i \in I_x$.
Therefore,the equivalence class
$\cartan_{(h,x)}^\ast$ is determined by the set $I_{(h,x)} = \{ i \in I_x \mid h \in H_i \}$.
Specifically,
\[
    \cartan_{(h,x)}^\ast = \left(\bigcup\limits_{i\in I_{(h,x)}} H_i\right) \cap
    \left(\bigcup\limits_{j\notin I_{(h,x)}} H_j\right)^c
\]
where $\,^c$ denote the complement (cf.~Subsection \ref{subsec-StratConnectedTorus} below).
The stratification of $\inertia{M}$ given by Equation \eqref{eq-StratTildeMDef}, then, is
given by sets of the form $\cartan_{(h,x)}^\ast\times M_{H_i} $
where $h \in H_i$ and $x \in M_{H_i}$. Intuitively, $\inertia{M}$ is
partitioned by isotropy types, and then further decomposed to
separate the closures of nearby strata with lower-dimensional fibers
in the $G$--direction.

As a particularly elucidating example, consider $G = \T^2 = \{ (s,
t) \mid s, t \in \T^1 \}$ and $M = \CP^2$ with action given by
\[
    (s, t)[z_0, z_1, z_2]
        =   [sz_0, tz_1, stz_2].
\]
Note that the points $[1, 0, 0]$, $[0, 1, 0]$, and $[0, 0, 1]$ are
fixed by $\T^2$. Near these three points, respectively, using
coordinates
\[
    (u_1, u_2) := \left(\frac{z_1}{z_0}, \frac{z_2}{z_0}\right), \;\;\;
    (v_0, v_2) := \left(\frac{z_0}{z_1}, \frac{z_2}{z_1}\right), \;\;\; \mbox{and}\;\;\;
    (w_0, w_1) := \left(\frac{z_0}{z_2}, \frac{z_1}{z_2}\right),
\]
the action is given by
\[
\begin{array}{rcl}
    (s, t)(u_1, u_2)    &=&     (s^{-1}t u_1, tu_2),                    \\
    (s, t)(v_0, v_2)    &=&     (st^{-1} v_0, s v_2), \;\;\mbox{and}    \\
    (s, t)(w_0, w_1)    &=&     (t^{-1} w_0, s^{-1} w_1).
\end{array}
\]
Note in particular that the action near each fixed point is
different, and hence the torus $\T^2$ is partitioned into $\sim$
classes in different ways at each.  However, the strata in
$\inertia{\CP^2}$ whose closures contain two fixed points have torus
fiber given by a subtorus of $\T^2$ whose partition into $\sim$
classes is compatible with both.

For instance, the torus fiber over $R = \{ [z_0, z_1, 0] \mid z_0,
z_1 \neq 0 \}$ in $\inertia{\CP^2}$ is the $1$--dimensional subtorus
$T$ of $\T^2$ consisting of points of the form $(s, s)$.  Any open
neighborhood of the orbit of a point in $R$ contains points with
trivial isotropy and points with isotropy $T$, so that $T \times R$
is partitioned into $\{ e \}\times R$ and $(T\smallsetminus \{ e
\})\times R$. It is easy to see that this partition is the
restriction to $T$ of the partitions of $\T^2$ at $[1,0,0]$ and
$[0,1,0]$.  Though it is not compatible with the partition at
$[0,0,1]$, this causes no difficulty as $[0,0,1]$ is separated from
the closure of $R$.

\subsubsection{The Adjoint Action}
\label{subsubsec-ExamplesAdjoint}

Let $G$ act on itself by conjugation. Then $\inertia{G} = \{ (h, k) \in G^2 \mid kh = hk \}$
is the set of commuting ordered pairs of elements of $G$ with diagonal $G$-action by conjugation.
The isotropy group of $(h, k)$ is given by $Z_G(h) \cap Z_G(k)$,
and the set of points fixed by $Z_G(h) \cap Z_G(k)$ is given by points of the form
$(l, j)$ where $l$ and $j$ are elements of the center of $Z_G(h) \cap Z_G(k)$.
Similarly, $\cartan_{(h,k)}^\ast$ consists of elements $j$ of $\cartan_{(h,k)}$
such that $Z_G(h)$ and $Z_G(j)$ coincide on a neighborhood of the $G$-orbit
$G(h,k)$.

\subsubsection{The Standard Action of $\SO{3}$ on $\R^3$}
\label{subsubsec-ExamplesSO3}

Let $G = \SO{3}$ act on $M = \R^3$ in the usual way.  For each point $x \in \R^3$ with
$x \neq 0$, we let $R_{x,\theta}$ with $\theta \in [0, 2\pi)$
denote rotation through the angle $\theta$ about the
line spanned by $x$ where we assume $\theta$ is a positive rotation with respect
to an oriented basis for $\R^3$ whose third element is $x$.  In particular,
$R_{x,0} = 1$ and $R_{x,\theta} = R_{-x, 2\pi - \theta}$  for each $x \in \R^3 \smallsetminus \{ 0\}$.
See \cite[Sec.~1.2 and 3.4]{duistermaatkolk} for a careful description of this
action, and note that our notation differs slightly to adapt to our situation.

There are three isotropy types that occur in $\inertia{\R^3}$: the
point $(e, 0)$ has isotropy group $\SO{3}$, points of the form $(R_{x,\pi}, 0)$
have conjugate isotropy groups isomorphic to $O(2)$, and all other points have
conjugate isotropy groups isomorphic to $\SO{2}$.  If $(h,x) \in \inertia{\R^3}$
such that $x \neq 0$ and $\cartan_{(h,x)} = G_{(h,x)} \cong \SO{2}$,
then any neighborhood of the orbit $G(h,x)$ small enough to not intersect
$\{ 0 \} \times \SO{3}$ contains only points with $\cartan_{(h,x)}$-isotropy
type $\cartan_{(h,x)}$ and $\{ e \}$.  Hence there are only two $\sim$ classes,
the identity and the nontrivial elements.
Similarly, if $(h,x) = (R_{x, \theta}, 0)$ with $\theta \neq \pi$,
it can be seen in a neighborhood of the orbit $G(h,x)$ small enough to contain no points of the
form $(R_{y,\pi}, 0)$ that $\cartan_{(h,x)}$ is as well partitioned into the same two
$\sim$ classes.  If $(h,x) = (R_{x, \pi}, 0)$, then as $G(h,x)$ contains $(R_{y,\pi}, 0)$
for each $y \in \R^3$ and as $R_{x,\pi}$ fixes $(R_{y,\pi}, 0)$ when $x$ and $y$ are orthogonal,
the torus $\cartan_{(h,x)}$ is partitioned into the $\sim$ classes $\{ e \}$, $\{ R_{x,\pi}\}$,
and $\{ R_{x, \theta} \mid \theta \in (0, \pi) \cup (\pi, 2\pi)\}$.  It follows that the maximal
decomposition of $\inertia{\R^3}$ induced by the stratification $\mathcal{S}_{(h,x)}$
consists of four sets:
\[
\begin{array}{rcl}
    S_1     &=&     \{ (e, 0) \},                                                       \\
    S_2     &=&     \{ (R_{x,\pi}, 0)\mid x \in \R^3\smallsetminus \{0\}\},             \\
    S_3     &=&     \{ (e, x)\mid x \in \R^3\smallsetminus \{0\} \},  \;\;\;\mbox{and}  \\
    S_4     &=&     \{ (R_{x,\theta}, x) \mid \theta \in (0, \pi) \cup (\pi, 2\pi), x \in \R^3 \}
                    \cup
                    \{ (R_{x, \pi}, x) \mid x \in \R^3 \smallsetminus \{ 0 \} \}.
\end{array}
\]
Note in particular that the map $\SO{3}\backslash\inertia{\R^3} \to \SO{3}\backslash\R^3$
given by $\SO{3}(h,x) \mapsto \SO{3}x$ is not a stratified mapping; for
$\theta \in (0,\pi)\cup(\pi,2\pi)$,
the points $\SO{3}(R_{x,\theta}, x)$ and $\SO{3}(R_{x,\theta}, 0)$ are mapped to points
with different isotropy types.

Similarly, consider the restriction of the $\SO{3}$-action to $M = \R^3\smallsetminus\{0\}$.
The maximal decomposition of $\inertia{M}$ given by Equation \eqref{eq-StratTildeMDef} has two
pieces,
\[
    \{ (e, x)\mid x \in \R^3\smallsetminus \{0\} \}
    \;\;\mbox{and}\;\; \{ (R_{x,\theta}, x)\mid \theta \in (0, 2\pi), x \in \R^3\smallsetminus \{0\} \}.
\]
Note in particular that in this case, $\inertia{M}$ is a smooth manifold with a single isotropy type,
and hence that this stratification is strictly finer than the stratification by isotropy types.

To understand this phenomenon, let $H$ be the subgroup of $\SO{3}$ isomorphic to $\SO{2}$
given by rotations about the $z$-axis.  Then considering the $H$-space $\R^3$
given by the restricted action, there are two isotropy types; points on the $z$-axis are fixed
by all of $H$, while points off the $z$-axis are fixed only by the identity.  It is easy to see,
then, that the partition of $\inertia{\R^3}$ given by the restriction of
the isotropy type stratification of $\R^3 \times H$ does not yield
a stratification of $\inertia{\R^3}$.  Hence, while the stratifications given by
Equations \eqref{eq-StratTildeMDef} and \eqref{eq-StratInertDef} are in general not the coarsest
stratifications of $\inertia{M}$ and $\inertia{X}$, they have the benefit
of giving a uniform, explicit stratification of the loop space $\inertia{M}$ and the
inertia space for all smooth $G$-manifolds under consideration.

%
%


\subsection{A Partition of Cartan Subgroups in Isotropy Groups}
\label{subsec-StratConnectedTorus}

In this subsection, we prove a  number of auxiliary results on topological properties of the
equivalence classes of the relation $\sim$ which has been defined in
Subsection \ref{subsec-ExamplesStratDescription}.
Throughout this section, let $Q$ be a smooth, not necessarily connected $G$-manifold and fix a closed
abelian subgroup $\cartan \leq G$ which need not be connected.
Assume that $Q$ is partitioned into a finite number of
$\cartan$-isotropy types.  We have in mind the case $Q = GV$
where $V$ is a slice for the $G$-action on $G \times M$ and $\cartan$
is a Cartan subgroup of the isotropy group of the origin in $V$, but we
state the results of this subsection more generally.

As above, for $s, t \in \cartan$, we say that $s \sim t$ when $Q^s = Q^t$.
Let $H_0, H_1, \ldots H_r$ be the finite collection of isotropy groups
for the action of $\cartan$ on $Q$. Then the $\sim$ class $[t]_{}$ of $t \in \cartan$ is given by
\[
    [t] = \bigcap\limits_{t \in H_i} H_i \cap  \left(\bigcup\limits_{t \not\in H_j} H_j\right)^c.
\]
That is, each $\sim$ class is determined by a subset of $\{ 1, 2, \ldots, r \}$; note that a
nonempty subset $I \subseteq \{ 1, 2, \ldots, r \}$ need not correspond to a nonempty
$\sim$ class.  Using this together with a dimension counting argument the following result
is derived immediately.

\begin{lemma}
\label{lem-StratConnectedTorusPartition}
The group $\cartan$ is partitioned into a finite number of $\sim$ classes.
Each $\sim$ class $[t]_{}$ is an open subset of the closed subgroup $t^\bullet$
of $\cartan$ defined by
\[
    t^\bullet
    :=
    \bigcap\limits_{t \in H_i} H_i
    =
    \bigcap\limits_{q \in Q^t} \cartan_q,
\]
and $\overline{[t]_{}}$ consists of a union of connected components of $t^\bullet$.
Moreover, each $\sim$ class has a finite number of connected components.
\end{lemma}

Also note the following.

\begin{lemma}
\label{lem-StratConnectedTorusFrontier}
Suppose $s, t \in \cartan$ such that $[s]_{} \cap \overline{[t]_{}} \neq \emptyset$.
Then for each connected component $[s]_{}^\circ$ of $[s]_{}$ and $[t]_{}^\circ$ of $[t]_{}$
such that $[s]_{}^\circ \cap \overline{[t]_{}^\circ} \neq \emptyset$
the relation $[s]_{}^\circ \subseteq \overline{[t]_{}^\circ}$ holds true.
\end{lemma}
\begin{proof}
Let $u \in [s]_{}^\circ \cap \overline{[t]_{}^\circ}$. Then $Q^s = Q^u$, and by continuity of the
action, $Q^t \subseteq Q^u$.  It follows that $Q^t \subseteq Q^s$, and hence that
$s^\bullet = \bigcap\limits_{q \in Q^s}
\cartan_q \leq \bigcap\limits_{q \in Q^t} \cartan_q = t^\bullet$.

Note that $[s]_{}^\circ$ is contained in a connected component
$(s^\bullet)^\circ$ of $s^\bullet$ which is contained in a connected component
$(t^\bullet)^\circ$ of $t^\bullet$.  Similarly, $\overline{[t]_{}}$ consists of
entire connected components of $t^\bullet$, so $u \in \overline{[t]_{}^\circ} = (t^\bullet)^\circ$.
Then $[s]_{}^\circ \subseteq (t^\bullet)^\circ = \overline{[t]_{}^\circ}$, completing the proof.
\end{proof}

For each $g \in G$, we let $\sim_g$ denote the equivalence relation defined
on $g\cartan g^{-1}$ in terms of its action on $Q$.  In particular,
if $g \in N_G(\cartan)$, then $\sim_g$ coincides with $\sim$.
It is easy to verify the following.

\begin{lemma}
\label{lem-StratConnectedTorusinG}
Let $s, t \in \cartan$.  Then $s \sim t$ if and only if $gsg^{-1} \sim_g gtg^{-1}$, i.e.
\[
    [gtg^{-1}]_g = g[t]_{}g^{-1},
\]
where $[-]_g$ denotes the equivalence class with respect to $\sim_g$.
\end{lemma}

Similarly, the following will be important when showing that certain $\sim$ classes
are sufficiently separated.

\begin{lemma}
\label{lem-StratConnectedSimsDisconnected}
Suppose $s, t \in \cartan$ such that $s \not\sim t$, and
$[s]_{}$ is diffeomorphic to $[t]_{}$.
Then $[s]_{} \cap \overline{[t]_{}} = \emptyset$.
\end{lemma}
\begin{proof}
If $Q^s \subseteq Q^t$,
then $\bigcap_{q \in Q^t} \cartan_q \leq \bigcap_{q \in Q^s} \cartan_q$,
so that $t^\bullet \leq s^\bullet$.  By Lemma \ref{lem-StratConnectedTorusPartition},
$[s]_{}$ and $[t]_{}$ are open subsets of $s^\bullet$ and $t^\bullet$, respectively,
so that as $[s]_{}$ and $[t]_{}$ are diffeomorphic,
$s^\bullet$ and $t^\bullet$ have the same dimension.  Additionally, $[s]_{}$
is open and dense in each connected component of $s^\bullet$ it intersects,
and similarly $[t]_{}$ in $t^\bullet$, so that $[s]_{}$ and
$[t]_{}$ do not intersect the same connected components of the closed
group $t^\bullet$.  The claim follows, and the argument is identical if
$Q^t \subseteq Q^s$.

So suppose $Q^s \not\subseteq Q^t$ and $Q^t \not\subseteq Q^s$.
If $l \in \overline{[s]_{}} \cap \overline{[t]_{}}$,
then by continuity of the action, $l$ fixes $Q^s \cup Q^t$.
Since $Q^t$ is a proper subset of $Q^s \cup Q^t$,
it follows that $l \not\sim s$.  Therefore,
$l \in \overline{[s]_{}} \smallsetminus [s]_{}$
and $[s]_{} \cap \overline{[t]_{}} = \emptyset$.
\end{proof}

For each $n \in N_G(\cartan)$,  conjugation by $n$
induces a diffeomorphism from $\cartan$ to itself which by Lemma \ref{lem-StratConnectedTorusinG}
acts on the set of $\sim$ classes. More precisely:

\begin{lemma}
\label{lem-StratConnectedNormalizerPermut}
The normalizer $N_G(\cartan)$ acts on the finite set of $\sim$ classes
in $\cartan$ in such a way that for each $n \in N_G(\cartan)$ and $t \in \cartan$,
the submanifold $n[t]_{} n^{-1}$ is diffeomorphic to $[t]_{}$.  Moreover, either
$n[t]_{} n^{-1} =  [t]_{}$ or $\overline{n[t]_{} n^{-1}} \cap [t]_{} = \emptyset$.
\end{lemma}


\subsection{Proof of Theorem \ref{thrm-StratConnected}}
\label{subsec-StratConnectedProof}

In this subsection, we prove Theorem
\ref{thrm-StratConnected}, establishing that for a $G$-manifold $M$ the germs $\mathcal{S}_{(h,x)}$ and
$\mathcal{R}_{G(h,x)}$ given by Equations \eqref{eq-StratTildeMDef} and \eqref{eq-StratInertDef}
define a smooth Whitney stratification of the loop space $\inertia{M}$ and a smooth stratification
of the inertia space $\inertia{X}$, respectively.

The general strategy is to first decompose
$\inertia{M}$ into its $G$-isotropy types.  Roughly speaking, isotropy types
consisting of smaller manifolds have larger $G$-fibers, so the fibers must further
be decomposed as illustrated in the examples in Subsection \ref{subsubsec-ExamplesAbelian}.  This is
accomplished by first decomposing a Cartan subgroup in the $G$-fiber into $\sim$ classes
using the results of Subsection \ref{subsec-StratConnectedTorus} and then partitioning nearby
by taking the $G$-orbits of these pieces.  A brief outline of the proof follows.

We begin with Lemma \ref{lem-finisotropytypes} which essentially guarantees that we can apply
the results of the preceding section on the $G$-saturation of a (linear) slice.
Then, in Lemma \ref{lem-StratConnectedSubgermIsotropyType} we confirm that the germ
$\mathcal{S}_{(h,x)}$ is that of a subset of $\inertia{M}$ consisting of points with
the same $G$-isotropy type.  Afterwards, we  prove
Lemmas \ref{lem-StratConnectedIndepTorus} and \ref{lem-StratConnectedSliceIndependent},
demonstrating that the germs $\mathcal{S}_{(h,x)}$ and hence the $\mathcal{R}_{G(h,x)}$
do not depend on the choices of the slice, the Cartan subgroup associated to $h$, and the
representative $(h,x)$ of the orbit $G(h,x)$.  With this,
we prove Proposition \ref{prop-StratConnectedManifolds}, showing that
$\mathcal{S}_{(h,x)}$ and $\mathcal{R}_{G(h,x)}$ are germs of smooth submanifolds of
$G \times M$.

With this, we are required to define a decomposition $\mathcal{Z}$ of a neighborhood
of $U$ of each point $(h,x) \in \inertia{M}$; indicating this decomposition and verifying
its properties involve the main technical details of the proof.
The definition of $\mathcal{Z}$ is given in Equation \eqref{eq-StratConnectedPieces} in a manner similar to
the stratification; the piece containing $(k, y)$ is defined in terms of the isotropy
type of $(k, y)$ and the $\sim$ class $\cartan_{(k,y)}^\ast$ of $k$ with respect to the
action near $G(k, y)$.  However, this definition
is given in terms of a slice at $(h,x)$ rather than $(k, y)$, so that
we must take into consideration the orbit of the $\sim$ class $\cartan_{(k,y)}^\ast$
under the action of the normalizer $N_{G_{(h,x)}}(\cartan_{(h,x)})$.
In particular, the pieces of $\mathcal{Z}$ are defined to be connected components so that,
though they are $G^\circ$-invariant, they need not be $G$-invariant.  However,
the $G$-action simply permutes the pieces of $\mathcal{Z}$ that are connected components
of the same $G$-invariant set.

As the definition of each piece of $\mathcal{Z}$ involves choosing a
particular point in each orbit near that of $(h,x)$ as well as a Cartan subgroup,
Lemmas \ref{lem-StratConnectedPieceWellDef} and \ref{lem-StratConnectedPieceIndepPoint}
demonstrate that the definition is independent of these choices and the resulting
partition is well-defined.  We then show in
Proposition \ref{prop-StratConnectedGermsCoincide} that the germs of the pieces
of the decomposition $\mathcal{Z}$ coincide with the stratification.  This in particular
requires a careful description of a $G$-invariant neighborhood $W$ of a $(k, y)$
small enough not to intersect certain $\sim$ classes in the Cartan subgroup $\cartan_{(k,y)}$.
Roughly speaking, $W$ is formed by removing the closures of the finite collection
of conjugates of $\cartan_{(k,y)}^\ast$ by $N_{G_{(h,x)}}(\cartan_{(k,y)}^\ast)$ from
the $G$-factor; it is on this neighborhood that the connected component of the
stratum containing $(k, y)$ coincides with the piece containing $(k, y)$.  As the stratum
containing $(k,y)$ has finitely many connected components in this neighborhood, it follows
that the germs coincide.  With this, we demonstrate that the partition
of a neighborhood of $(h,x)$ is finite in Lemma \ref{lem-StratConnectedLocFinite},
that it satisfies the condition of frontier in Proposition \ref{prop-StratConnectedFrontier},
and that it satisfies Whitney's condition B in Proposition \ref{prop-StratConnectedWhitneyB}.
This completes the outline, and we now proceed with the proof.

First we assume to have fixed a $G$-invariant riemannian metric on $M$, a bi-invariant metric
on $G$, and that $G \times M$ carries the product metric. By $(h,x)$ we will always denote a
point of the loop space $\inertia{M}$, and by $V_{(h,x)}$ a linear slice in $G\times M$ at  $(h,x)$.
The isotropy group $G_{(h,x)} = Z_{G_x} (h)$ of $(h,x)$ will be denoted by $H$, and the normal space
$T_{(h,x)} (G\times M) / T_{(h,x)} (G(h,x))$ by $N_{(h,x)}$.

\begin{lemma}
\label{lem-finisotropytypes}
  Let $K$ be a closed subgroup of $G$, and $V_{(h,x)}$ a (linear) slice for the $G$-action on
  $G \times M$ as above. Then the $K$-manifold $Q:=GV_{(h,x)}$ has a finite number of $K$-isotropy types.
\end{lemma}
\begin{proof}
Let $\Psi: V_{(h,x)} \rightarrow N_{(h,x)}$ denote an $H$-invariant embedding of the slice $V_{(h,x)}$
into the normal space $N_{(h,x)}$ such that its image is an open convex neighborhood of the origin.
Choose an  $H$-invariant open convex neighborhood $B$ of the origin of $N_{(h,x)}$
which is relatively compact in $\Psi \big( V_{(h,x)} \big) $.
For each point $(k, y) \in GV_{(h,x)}$ choose a slice $Y_{(k,y)}$ for the $K$-action on $GV_{(h,x)}$.
Then the family $\{ KY_{(k,y)} \}_{(k,y) \in GV_{(h,x)}}$
is an open cover of $G\Psi^{-1} ( \overline{B})$ which has to admit a finite subcover by compactness
of $G\Psi^{-1} ( \overline{B})$.  Since each $KY_{(k,y)}$ has a finite number of $K$-isotropy types
by \cite[Lem.~4.3.6]{PflaumBook}, it follows that $G\Psi^{-1} ( \overline{B})$, hence $G\Psi^{-1} (B)$
has a finite number of $K$-isotropy types. However, $GV_{(h,x)}$ contains the same isotropy types
as $G\Psi^{-1} (B)$, since the action by $t\in (0,1]$ on  $V_{(h,x)}$ is $G$-equivariant, and for each
$v \in V_{(h,x)}$  there is  a $t\in (0,1]$ with $tv \in \Psi^{-1} (B)$.
Hence $GV_{(h,x)}$ itself has a finite number of $K$-isotropy types.
\end{proof}

The Lemma implies in particular that the results of Subsection \ref{subsec-StratConnectedTorus} apply
to $Q = GV_{(h,x)}$ for each abelian subgroup $\cartan = K$ of $G$.

\begin{lemma}
\label{lem-StratConnectedSubgermIsotropyType}
The set germ $\mathcal{S}_{(h,x)}$ is contained in the set germ at $(h,x)$ of points of $\inertia{M}$
having the same isotropy type as $(h,x)$ with respect to the $G$-action on $G \times M$.
\end{lemma}
\begin{proof}
Suppose
\[
    (k, y) \in V_{(h,x)}^H \cap (\cartan_{(h,x)}^\ast \times M).
\]
Since $H$ fixes $(k,y)$ and $k \in \cartan_{(h,x)} \leq H$, one obtains $ky = y$ and
$(k, y) \in \inertia{M}$. By $G$-invariance of $\inertia{M}$ we get $G(k,y) \subseteq \inertia{M}$,
hence $\mathcal{S}_{(h,x)}$ is the germ of a subset of $\inertia{M}$.  Now observe that
$V_{(h,x)}^H = (V_{(h,x)})_H \subseteq (G \times M)_H$,
where $(V_{(h,x)})_H$ and $(G \times M)_H$  denote the subsets of points having isotropy group $H$.
Hence the isotropy group of every point in the $G$-orbits defining
$\mathcal{S}_{(h,x)}$ is conjugate to $H$, and
$\mathcal{S}_{(h,x)}$ is a subgerm of $(G \times M)_{(H)}$.
\end{proof}

The following two lemmas demonstrate that the stratification $\mathcal{S}_{(h,x)}$
does not depend on the choice of a Cartan subgroup $\cartan_{(h,x)}$ nor on the particular choice
of a slice $V_{(h,x)}$.

\begin{lemma}
\label{lem-StratConnectedIndepTorus}
Let $(h,x) \in \inertia{M}$ with isotropy group $H = Z_{G_x}(h)$.  The germ
$\mathcal{S}_{(h,x)}$ does not depend on the choice of the Cartan subgroup
$\cartan_{(h,x)}$ of $H$.
\end{lemma}
\begin{proof}
Suppose $\cartan_{(h,x)}$ and $\cartan_{(h,x)}^\prime$ are two Cartan subgroups
of $H$ associated to $h$.  Then $\cartan_{(h,x)}/\cartan_{(h,x)}^\circ$ is generated
by $h\cartan_{(h,x)}^\circ$ and $\cartan_{(h,x)}^\prime/\cartan_{(h,x)}^{\prime\,\circ}$
is generated by $h\cartan_{(h,x)}^{\prime\,\circ}$, where here and in the rest of this section $K^\circ$
denotes the connected component of the neutral element in a Lie group $K$.
Under the correspondence given in \eqref{eq-StratConnectedCartanConjugate} both
$\cartan_{(h,x)}/\cartan_{(h,x)}^\circ$ and $\cartan_{(h,x)}^\prime/\cartan_{(h,x)}^{\prime\,\circ}$ then
correspond to $\langle hH^\circ\rangle \leq H/H^\circ$. It follows by
\cite[IV.~Prop.~4.6]{BroeckertomDieck} that $\cartan_{(h,x)}$ and $\cartan_{(h,x)}^\prime$ are conjugate,
so that there is a $g \in H$ such that $g\cartan_{(h,x)}g^{-1} = \cartan_{(h,x)}^\prime$.  Then,
as $g \in H$, the space $V_{(h,x)}^H$ is left invariant by $g$, hence $ghg^{-1} = h$.  Therefore, if
$k \in \cartan_{(h,x)}$ with $k \sim h$ as elements of $\cartan_{(h,x)}$
acting on $GV_{(h,x)}$, Lemma \ref{lem-StratConnectedTorusinG} implies that
$gkg^{-1} \sim ghg^{-1} = h$ as elements of $\cartan_{(h,x)}^\prime$
acting on $GV_{(h,x)}$.  It follows that conjugation by $g$ induces a
diffeomorphism of $\cartan_{(h,x)}^\ast$ onto $\cartan_{(h,x)}^{\prime\,\ast}$,
so that
\[
    g\left(V_{(h,x)}^H \cap \left(\cartan_{(h,x)}^\ast\times M \right)\right)
    =
    \left(V_{(h,x)}^H \cap \left(\cartan_{(h,x)}^{\prime\,\ast}\times M \right)\right),
\]
and
\[
    G\left(V_{(h,x)}^H \cap \left( \cartan_{(h,x)}^\ast \times M\right)\right)
    =
    G\left(V_{(h,x)}^H \cap \left( \cartan_{(h,x)}^{\prime\,\ast} \times M\right)\right).
\]
\end{proof}

\begin{lemma}
\label{lem-StratConnectedSliceIndependent}
The germ $\mathcal{S}_{(h,x)}$ is independent of the particular choice of the slice $V_{(h,x)}$ at $(h,x)$.
\end{lemma}
\begin{proof}
Suppose $V_{(h,x)}$ and $W_{(h,x)}$ are two choices of linear slices at $(h,x)$ for the $G$-action on $G \times M$.
Let $\cartan_{(h,x)}$ be Cartan subgroup of $H$ associated to $h$. By
Lemma \ref{lem-StratConnectedIndepTorus}, we may assume that the stratum containing $(h,x)$ is defined
with respect to each of the two slices using this Cartan subgroup.
Note that by the slice theorem and the assumptions on $V_{(h,x)}$ and $W_{(h,x)}$ the open sets
$GV_{(h,x)} \cong G \times_H V_{(h,x)}$ and $GW_{(h,x)} \cong G \times_H W_{(h,x)}$ are $G$-diffeomorphic
and hence $\cartan_{(h,x)}$-diffeomorphic.  Therefore,
the $\sim$ classes in $\cartan_{(h,x)}$ do not depend on the choice of the slice.

Letting
$\mathcal{N} = \{n \in N_G(\cartan_{(h,x)}) : n\cartan_{(h,x)}^\ast n^{-1} \neq \cartan_{(h,x)}^\ast\}$,
we have by Lemma \ref{lem-StratConnectedNormalizerPermut} that the set
\[
    C = (H \smallsetminus hH^\circ) \cup
    \bigcup\limits_{n \in \mathcal{N}}
    n\left(\overline{\cartan_{(h,x)}^\ast} \right) n^{-1}
\]
is a closed subset of $G$ disjoint from $\cartan_{(h,x)}^\ast$. Hence
$V_{(h,x)} \cap (C \times M)^c$ is an open neighborhood of $(h,x)$ in $V_{(h,x)}$.
We may therefore assume after possibly shrinking $V_{(h,x)}$ and $W_{(h,x)}$ that
$V_{(h,x)} \cap (C \times M) = W_{(h,x)} \cap (C \times M) = \emptyset$.
Clearly, shrinking the slice does not affect the germ of the
stratum.  With this, we let $O$ denote the $G$-invariant open neighborhood
$O := G V_{(h,x)} \cap G W_{(h,x)}$
of $(h,x)$ and claim that
\[
    O \cap G \left( V_{(h,x)}^H \cap \left(\cartan_{(h,x)}^\ast \times M\right) \right)
    =
    O \cap G \left( W_{(h,x)}^H \cap \left(\cartan_{(h,x)}^\ast\times M\right) \right).
\]

Any element of
$O \cap G \left( V_{(h,x)}^H \cap \left(\cartan_{(h,x)}^\ast \times M\right) \right)$
is in the $G$-orbit of some
$(k, y)\in O \cap V_{(h,x)}^H \cap \left(\cartan_{(h,x)}^\ast \times M\right)$.
As $(k, y) \in O$, there is a $g \in G$ such that $g(k,y) \in W_{(h,x)}$.
Since $G_{(k, y)} = H$ and $G_{g(k,y)} = gG_{(k,y)}g^{-1} \leq H$,
\cite[Lem.~4.2.9]{PflaumBook} implies that $G_{g(k,y)} = H$.
In particular, $g \in N_G(H)$, and $g(k,y) \in W_{(h,x)}^H$.

Now, $k \in \cartan_{(h,x)}^\ast \subseteq \cartan_{(h,x)}$ by definition,
so that $gkg^{-1}\in g\cartan_{(h,x)}g^{-1}$.  As $g \in N_G(H)$, it follows that
$g\cartan_{(h,x)}g^{-1} \leq H$.  Noting that
$k$ is an element of the connected set $\cartan_{(h,x)}^\ast \ni h$,
we have that $k$ is in the same connected component of $\cartan_{(h,x)}$ as $h$
and so $\cartan_{(h,x)}$ is a Cartan subgroup of $H$ associated to $k$ as well as $h$.
It is then easy to see that $g\cartan_{(h,x)}g^{-1}$
is a Cartan subgroup of $H$ associated to $gkg^{-1}$.
Moreover, because $W_{(h,x)}$ is disjoint from $(H \smallsetminus hH^\circ)\times M \subseteq C \times M$,
it must be that $gkg^{-1} \in hH^\circ$.  By
\cite[IV.~Prop.~4.6]{BroeckertomDieck}, there is a $\tilde{h} \in H$
such that $\tilde{h}g\cartan_{(h,x)}g^{-1}\tilde{h}^{-1} = \cartan_{(h,x)}$, and hence
$\tilde{h}g\in N_G(\cartan_{(h,x)})$.  That $\tilde{h} \in H = G_{g(k,y)}$ implies
$\tilde{h}gkg^{-1}\tilde{h}^{-1} = gkg^{-1}$, so that
$gkg^{-1}    \in     \tilde{h}g\cartan_{(h,x)}g^{-1}\tilde{h}^{-1} = \cartan_{(h,x)}$.
Moreover, as $k \in \cartan_{(h,x)}^\ast$, we have in addition that
$gkg^{-1} = \tilde{h}gkg^{-1}\tilde{h}^{-1} \in \tilde{h}g\cartan_{(h,x)}^\ast g^{-1}\tilde{h}^{-1}$.
Therefore,
\[
    g(k,y) \in O \cap W_{(h,x)}^H \cap \left(\tilde{h}g \cartan_{(h,x)}^\ast g^{-1}\tilde{h}^{-1}
    \times M\right).
\]
However, as $W_{(h,x)} \cap (C \times M) = \emptyset$, as $C$ contains all of the nontrivial
conjugates of elements of $\cartan_{(h,x)}^\ast$ by elements of $N_G(\cartan_{(h,x)})$,
and as $\tilde{h}g \in N_G(\cartan_{(h,x)})$, it must be that
$\tilde{h}g \cartan_{(h,x)}^\ast g^{-1}\tilde{h}^{-1} = \cartan_{(h,x)}^\ast$.  Hence,
\[
    g(k,y) \in O \cap W_{(h,x)}^H \cap \left( \cartan_{(h,x)}^\ast \times M\right).
\]
Switching the roles of $W_{(h,x)}$ and $V_{(h,x)}$ completes the proof.
\end{proof}

Note that if $(h,x) \in \inertia{M}$ and $g \in G$, then $gV_{(h,x)}$ is a slice
at $g(h,x)$, $g(V_{(h,x)}^H) = (gV_{(h,x)})^{gHg^{-1}}$, and $g\cartan_{(h,x)}g^{-1}$
is a Cartan subgroup of $gHg^{-1}$ associated to $ghg^{-1}$.  Therefore, as
the $\sim$ classes depend only on the action on $GV_{(h,x)}$, Lemmas
\ref{lem-StratConnectedIndepTorus} and \ref{lem-StratConnectedSliceIndependent}
imply that $g\mathcal{S}_{(h,x)} = \mathcal{S}_{g(h,x)}$, so that in particular
$\mathcal{R}_{G(h,x)}$ is well-defined. \vspace{2mm}

Now we have the means to verify the following crucial result.

\begin{proposition}
\label{prop-StratConnectedManifolds}
Each $\mathcal{S}_{(h,x)}$ is the germ of a smooth $G$-submanifold of $G \times M$,
and each $\mathcal{R}_{(h,x)}$ is the germ of a  smooth submanifold of
$G\backslash(G \times M)$.
\end{proposition}
Of course, $G\backslash(G \times M)$ is not itself a smooth manifold but rather a differentiable
space.  By a smooth submanifold of $G\backslash(G \times M)$, we mean a differentiable subspace
of $G\backslash(G \times M)$ that is itself a smooth manifold.
\begin{proof}
Since the germ $\mathcal{S}_{(h,x)}$ does not depend on the choice of a particular slice by Lemma
\ref{lem-StratConnectedSliceIndependent}, we choose the slice $V_{(h,x)}$ at $(h,x)$
to be the image under the exponential map of an open ball $B_{(h,x)}$ around the origin of the
normal space $N_{(h,x)}$. Note that $N_{(h,x)}$ naturally carries an $H$-invariant inner product since
we have
initially  fixed an invariant riemannian metric on $M$ and a bi-invariant riemannian metric on $G$.

Since $(G \times_H V_{(h,x)})^H$ is a totally geodesic submanifold of $G \times_H V_{(h,x)}$
by \cite[6.1]{michor}, the exponential map at $(h,x)$ maps
$B_{(h,x)}^H =  N_{(h,x)}^H \cap B_{(h,x)}$ onto $V_{(h,x)}^H$.  Similarly, as $\cartan_{(h,x)}^\ast$
is an open subset of the closed subgroup
$h^\bullet$ of $H$ by Lemma \ref{lem-StratConnectedTorusPartition}, the relation
$T_h \cartan_{(h,x)}^\ast = T_h h^\bullet$ holds true.
It follows that the exponential map associated to the product metric on $G \times M$
maps the subspace
\[
    N_{(h,x)}^H \cap (T_h h^\bullet\oplus T_x M) \cap B_{(h,x)}
\]
onto $V_{(h,x)}^H \cap \left(\cartan_{(h,x)}^\ast \times M\right)$, which
is then diffeomorphic to an open neighborhood of the origin in a linear space.

Noting that $G \times_H V_{(h,x)}^H \cong G/H \times V_{(h,x)}^H$, the $G$-diffeomorphism
\[
    \Psi    \co     G \times_H V_{(h,x)}
                    \longrightarrow     GV_{(h,x)}
                    \subseteq           G \times M
\]
induced by the exponential map restricts to a $G$-diffeomorphism
\[
    G/H \times \left( V_{(h,x)}^H \cap \left(\cartan_{(h,x)}^\ast\times M\right)\right)
    \longrightarrow
    G\left( V_{(h,x)}^H \cap \left(\cartan_{(h,x)}^\ast\times M\right)\right).
\]
Moreover, $\Psi$ induces a map on quotient spaces which is a homeomorphism
\[
    G\backslash\left( G/H \times
    \left( V_{(h,x)}^H \cap \left(\cartan_{(h,x)}^\ast\times M\right)\right)\right)
    \longrightarrow
    G\backslash\left(G\left( V_{(h,x)}^H \cap
    (\cartan_{(h,x)}^\ast\times M)\right)\right).
\]
Hence,
\[
    V_{(h,x)}^H \cap \left(\cartan_{(h,x)}^\ast\times M\right)
    \cong
    G\backslash\left(G\left( V_{(h,x)}^H \cap
    (\cartan_{(h,x)}^\ast\times M)\right)\right)
\]
is a topological submanifold of $G\backslash (G \times M)$.
On the differentiable space $G\backslash (G \times M)$, the structure sheaf
$\mathcal{O}^\infty_{G\backslash(G \times M)}$ is locally that of $G$-invariant functions on
$G \times M$ (see Section \ref{sec-StrucSheaf}).  Similarly, the $G$-invariant
$\mathcal{C}^\infty$ functions on
$G/H \times \left( V_{(h,x)}^H \cap \left(\cartan_{(h,x)}^\ast\times M\right)\right)$
are exactly the $\mathcal{C}^\infty$ functions on
$V_{(h,x)}^H \cap \left(\cartan_{(h,x)}^\ast\times M\right)$ by
\cite[Prop.~5.2]{tomdieck} and \cite[Thm.~1.22 (5)]{NGonzalezSanchoBook}.
Therefore, $G\backslash\left(G\left( V_{(h,x)}^H \cap (\cartan_{(h,x)}^\ast\times M)\right)\right)$,
whose set germ at $(h,x)$ coincides with $\mathcal{R}_{(h,x)}$, is a smooth
submanifold of the differentiable space $G\backslash(G \times M)$.
\end{proof}

In order for the germs $\mathcal{S}_{(h,x)}$ to define a stratification,
one must verify that for each $(h,x) \in \inertia{M}$ there is a neighborhood
$U$ in $\inertia{M}$ and a decomposition $\mathcal{Z}$ of $U$ such that for
all $(k, y) \in \inertia{M}$, the germ $\mathcal{S}_{(k,y)}$
coincides with the germ of the piece of $\mathcal{Z}$ containing
$(k, y)$.  Set $U := GV_{(h,x)} \cap \inertia{M}$. We now define the decomposition
$\mathcal{Z}$ of $U$.
Given  $(\tilde{k}, \tilde{y}) \in U$ there is a $\tilde{g} \in G$ such that
$\tilde{g}(\tilde{k}, \tilde{y}) \in V_{(h,x)}$. Put $(k, y) = \tilde{g}(\tilde{k}, \tilde{y})$
and $K = G_{(k,y)} \leq H$, and let $\cartan_{(k,y)}$ be a Cartan subgroup in $K$
associated to $k$.  We define the piece of $\mathcal{Z}$ containing $(\tilde{k}, \tilde{y})$
to be the connected component containing $(\tilde{k}, \tilde{y})$ of the set
$\mathcal{U}_{\tilde{g}}^{\cartan_{(k,y)}}(\tilde{k}, \tilde{y})$ which is defined as the $G$-saturation
of the set of points $(l, z) \in (V_{(h,x)})_K \cap (\cartan_{(k,y)} \times M)$ such that
$\cartan_{(k,y)}$ is a Cartan subgroup of $K$ associated to $l$
and such that the $\sim$ class $\cartan_{(l, z)}^\ast$ is conjugate to $\cartan_{(k,y)}^\ast$
via an element of $N_H (\cartan_{(k,y)})$.  Note that as above $\cartan_{(l, z)}^\ast$ is
the $\sim$ class of $l$ in $\cartan_{(l, z)} = \cartan_{(k,y)}$
with respect to its action on $GV_{(l, z)}$, where $V_{(l, z)}$ is a slice for the $G$-action of
on $G \times M$ at $(l, z)$.  Observe that by \cite[Proposition 1.3(2)]{Schwarz}.
as the action of $K$ on $V_{(h,x)}$ is linear, the slice representation of each point
in $V_{(h,x)}$ with isotropy group $K$ is isomorphic.  It follows that the set
$\mathcal{U}_{\tilde{g}}^{\cartan_{(k,y)}}(\tilde{k}, \tilde{y})$ can be written as
\begin{equation}
  \label{eq:decUrep}
  \mathcal{U}_{\tilde{g}}^{\cartan_{(k,y)}}(\tilde{k}, \tilde{y}) =
  G \Big( \bigcup\limits_{n \in N_H (\cartan_{(k,y)})} (V_{(h,x)})_K \cap \big( n \cartan_{(k,y)}^\ast n^{-1} \times M \big) \Big)
\end{equation}

Observe that $(V_{(h,x)})_K$ is closed under the action of scalars $\in (0,1]$  and open in $V_{(h,x)}^K$
as a consequence of Lemma \ref{lem-finisotropytypes}. Moreover,
\[
    (V_{(h,x)})_K = V_{(h,x)}^K \setminus \bigcup\limits_{\{ (H_v)  \mid v\in V_{(h,x)} \& K < H_v \}}
    (V_{(h,x)})_{(H_v)} ,
\]
where the union is over the (finitely many) isotropy classes that (properly) contain $K$.
An analogous relation holds true for $(N_{(h,x)})_K$.
Since all fixed point sets $N_{(h,x)}^{H_v}$ are algebraic, the set $(N_{(h,x)})_K \cap S_{(h,x)}$
is a semialgebraic subset of $N_{(h,x)}$, where  $S_{(h,x)}$ is a sphere in $N_{(h,x)}$.
Hence, $(N_{(h,x)})_K \cap S_{(h,x)}$ and $(N_{(h,x)})_K$ have finitely many connected components by
\cite[Sec.~2.4]{BocCosRoyRAG}. Therefore, $(V_{(h,x)})_K$ has finitely many components, too.

Next, we want to show that $(V_{(h,x)})_K \cap \big( n \cartan_{(k,y)}^\ast n^{-1} \times M \big)$
for $n \in N_H (\cartan_{(k,y)}) $ has finitely many components as well.
To this end note first that for each $(\tilde{k}, \tilde{y}) \in V_{(h,x)}$ the group element
$\tilde{k}$ lies in the same connected component $hH^\circ$ of $H$ as $h$, since
$V_{(h,x)}$ is connected. Therefore, for any Cartan subgroup $\cartan_{(\tilde{k},\tilde{y})}$ of the isotropy
group $K$ of $(\tilde{k}, \tilde{y})$ associated to $\tilde{k}$, $\cartan_{(\tilde{k},\tilde{y})}$ is conjugate
to a subgroup of $\cartan_{(h,x)}$ in $H$.  To see this, note that
by \cite[IV.~Prop.~4.2]{BroeckertomDieck} and its proof, a Cartan subgroup
of $K$ associated to $\tilde{k}$ is generated by a maximal torus in $Z_K(\tilde{k})$ and $\tilde{k}$,
while a Cartan subgroup of $H$ associated to $\tilde{k}$ is generated by a maximal torus
in $Z_H(\tilde{k})$ and $\tilde{k}$.  Say
$\tilde{h}\cartan_{(\tilde{k},\tilde{y})}\tilde{h}^{-1} \leq \cartan_{(h,x)}$
for some $\tilde{h} \in H$.
Then $(k,y) := \tilde{h}(\tilde{k},\tilde{y}) \in V_{(h,x)}$ and
$k \in \cartan_{(h,x)}$.  It follows that we can
always choose a representative $(k,y)\in V_{(h,x)}$ from an orbit
such that $\cartan_{(k,y)} \leq \cartan_{(h,x)}$.

Since $(V_{(h,x)})_K$ is closed under multiplication by
scalars $t\in (0, 1]$ using the linear structure it inherits from $N_{(h,x)}$,
each point $t(k, y)$  has isotropy group $K$, hence its
$G$-coordinate lies in $K$.  As $K$ and hence $K\times M$
is closed, it then must be that $\lim_{t\to 0} t(k,y) = (h,x) \in K \times M$,
which means $h \in K$.  In particular, $h$ and $k$ are in
the same connected component of $K$ and hence that the Cartan subgroup
$\cartan_{(k,y)}$ of $K$ associated to $k$ is conjugate in $K$ to a Cartan
subgroup of $K$ associated to $h$.  It follows in particular that there
is a $\tilde{k} \in K$ such that $h \in \tilde{k}\cartan_{(k,y)}\tilde{k}^{-1}$.
However, this implies that $\tilde{k}^{-1}h\tilde{k} \in \cartan_{(k,y)}$,
so that as $\tilde{k} \in K \leq H$ fixes $h$, we have from the beginning that
$h \in \cartan_{(k,y)}$.

Now, recall that the equivalence class $\cartan_{(h,x)}^\ast$ is an
open subset of the closed subgroup $h^\bullet$ of $\cartan_{(h,x)}$. With respect
to the action of $\cartan_{(h,x)}$ on $GV_{(h,x)}$, there are a
finite number of $\sim$ classes in $\cartan_{(h,x)}$.  Hence, there is a neighborhood
$O$ of $h$ in $G$ which only intersects $\sim$ classes in
$\cartan_{(h,x)}$ whose closures contain $h$.  Assume $V_{(h,x)}$ is
a slice chosen small enough so that $V_{(h,x)} \subseteq O \times
M$, and pick $(k, y) \in V_{(h,x)}$. Then, one can chose the Cartan subgroup
$\cartan_{(k,y)} \leq G_{(k, y)}$ with $k \in \cartan_{(k,y)} $ such that
$\cartan_{(k,y)} \leq \cartan_{(h,x)}$.  Since the slice $V_{(k,y)}$ at $(k,y)$ may be
shrunk such that $GV_{(k,y)} \subseteq GV_{(h,x)}$, it follows from the definition of
$\sim$ that $\cartan_{(k,y)}^\ast$ is the intersection of a union of $\sim$ classes in
$\cartan_{(h,x)}$ with $\cartan_{(k,y)}$. In particular, as the closure of each such $\sim$
classes contains $h$, $\overline{\cartan_{(k,y)}^\ast}$ and $\cartan_{(k,y)}$ both contain $h$.
By the proof of Proposition \ref{prop-StratConnectedManifolds},
the relation $T_h (n \overline{\cartan_{(k,y)}^\ast} n^{-1}) = T_h (n k^\bullet n^{-1})$ holds true
for all $n \in  N_H (\cartan_{(k,y)})$, where $k^\bullet$ is the intersection of all
isotropy groups of the $\cartan_{(k,y)}$-action on $GV_{(k,y)}$ which contain $k$.
It follows that the exponential map associated to the
product metric on $G \times M$ maps the subspace
\begin{equation}
  \label{eq:modspace}
  (N_{(h,x)})_K \cap \big( T_h (n k^\bullet n^{-1}) \oplus T_x M \big) \cap B_{(h,x)}
\end{equation}
onto $(V_{(h,x)})_K \cap \left(n \overline{\cartan_{(k,y)}^\ast} n^{-1} \times M\right)$.
By construction, \eqref{eq:modspace} is a semialgebraic subset of $N_{(h,x)}$, and invariant under the action
of $t\in (0,1]$.

Let us now describe the preimage of
$(n \cartan_{(k,y)}^\ast n^{-1} \times M)$ under the exponential map. Since there are only
finitely many $\sim$ classes in $\cartan_{(k,y)}$, one can find finitely many elements
$l_1,\ldots ,l_\alpha \in \cartan_{(k,y)}$
such that each group $l_\iota^\bullet$, $\iota=1,\ldots,\alpha$  has dimension less than $\dim k^\bullet$,
and such that
\[
   \cartan_{(k,y)}^\ast = \overline{\cartan_{(k,y)}^\ast} \setminus \bigcup_{\iota =1}^\alpha l_\iota^\bullet \: .
\]
This implies that $\exp$ maps the set
\begin{equation}
  \label{eq:modspace2}
  (N_{(h,x)})_K \cap \Big( \big( T_h (n k^\bullet n^{-1}) \setminus \bigcup_{\iota =1}^\alpha
  T_h (n l_\iota^\bullet n^{-1}) \big) \oplus T_x M \Big) \cap B_{(h,x)}
\end{equation}
onto $(V_{(h,x)})_K \cap \big( n \cartan_{(k,y)}^\ast n^{-1} \times M \big)$.
But \eqref{eq:modspace2} is semialgebraic by construction, and invariant under the action of $t\in (0,1]$.
Hence, \eqref{eq:modspace2} and thus $(V_{(h,x)})_K \cap \big( n \cartan_{(k,y)}^\ast n^{-1} \times M \big)$
have both finitely many connected components, and are invariant under the action of $t\in (0,1]$, too.
Since $G$ is compact, and since there are only finitely many different sets $n \cartan_{(k,y)}^\ast n^{-1} $,
when $n$ runs through the elements of $N_H (\cartan_{(k,y)}) $, Eq.~\eqref{eq:decUrep} entails the following.

\begin{lemma}
\label{lem-StratConnectedUkyFiniteCC}
Suppose $V_{(h,x)}$ is given by the image under the exponential map of a sufficiently small
ball $B_{(h,x)}$ in the normal space $N_{(h,x)}$, and $(k,y)\in V_{(h,x)}$. Then the set
$\exp_{(h,x)}^{-1}\big(\mathcal{U}_{\tilde{g}}^{\cartan_{(k,y)}}\big( G(k,y) \big)\big) \cap B_{(h,x)}$
is invariant under multiplication by scalars in $(0, 1]$.
Moreover, each set $\mathcal{U}_{\tilde{g}}^{\cartan_{(k,y)}}\big( G(k,y) \big)$ has a finite number
of connected components.
\end{lemma}

At the moment, the set $\mathcal{U}_{\tilde{g}}^{\cartan_{(k,y)}}(\tilde{k}, \tilde{y})$
appears to depend both on the choice of $(k, y) \in V_{(h,x)}$ in the $G$-orbit of
$(\tilde{k}, \tilde{y})$ and the Cartan subgroup $\cartan_{(k, y)}$.  With the following two lemmas,
we will demonstrate that this is not the case, allowing us to simplify the notation.

\begin{lemma}
\label{lem-StratConnectedPieceWellDef}
The set $\mathcal{U}_{\tilde{g}}^{\cartan_{(k,y)}}(\tilde{k}, \tilde{y})$ does not depend
on the particular representative $(k, y) \in V_{(h,x)}$ of the $G$-orbit of $(\tilde{k}, \tilde{y})$,
hence does not depend on $\tilde{g}$.
\end{lemma}
\begin{proof}
Suppose $g \in G$ also satisfies $g(\tilde{k}, \tilde{y}) =: (k^\prime, y^\prime) \in V_{(h,x)}$.
It follows that $(k^\prime, y^\prime) = g\tilde{g}^{-1}(k,y)$,
so that we have $g\tilde{g}^{-1} \in H$ by (SL\ref{It:SlInv}), since
$g\tilde{g}^{-1}V_{(h,x)} \cap V_{(h,x)} \neq \emptyset$.
Let $\tilde{h} = g\tilde{g}^{-1}$. Then the isotropy group of $(k^\prime, y^\prime)$
is $K^\prime := \tilde{h}K\tilde{h}^{-1} \leq H$.
Let $\cartan_{(k^\prime, y^\prime)}$ be a choice of a Cartan subgroup
in $K^\prime$ associated to $k^\prime$.  Note that since $\tilde{h}\cartan_{(k,y)}\tilde{h}^{-1}$
is clearly a Cartan subgroup of $K^\prime$ associated to $k^\prime$ as well,
there exists by \cite[IV.~Prop.~4.6]{BroeckertomDieck}
a $\tilde{k} \in K^\prime$ such that
$\cartan_{(k^\prime, y^\prime)} = \tilde{k}\tilde{h}\cartan_{(k, y)} \tilde{h}^{-1}\tilde{k}^{-1}$.
Moreover, by Lemma \ref{lem-StratConnectedTorusinG},
$\cartan_{(k^\prime, y^\prime)}^\ast =
\tilde{k}\tilde{h}\cartan_{(k, y)}^\ast \tilde{h}^{-1}\tilde{k}^{-1}$.

Let $(\tilde{l}, \tilde{z}) \in \mathcal{U}_{\tilde{g}}^{\cartan_{(k,y)}}(\tilde{k}, \tilde{y})$.
Then there exists a $\tilde{g}^\prime \in G$ such that
$(l, z) := \tilde{g}^\prime(\tilde{l}, \tilde{z}) \in (V_{(h,x)})_K$
and an $n \in N_H(\cartan_{(k,y)})$ such that
 $\cartan_{(l, z)}^\ast  = n \cartan_{(k,y)}^\ast n^{-1}$.
As $(l, z)$ has isotropy group $K$, one has
$\tilde{h}(l, z) \in (V_{(h,x)})_{K^\prime}$.
Since $\tilde{k} \in K^\prime$, we get $\tilde{k}\tilde{h}(l, z) = \tilde{h}(l, z)$.
Again by Lemma \ref{lem-StratConnectedTorusinG}, we therefore obtain
\[
\begin{array}{rcl}
        \cartan_{\tilde{h}(l, z)}^\ast
        &=&
        \cartan_{\tilde{k}\tilde{h}(l, z)}^\ast            \\
        &=&
        \tilde{k}\tilde{h}\cartan_{(l, z)}^\ast \tilde{h}^{-1}\tilde{k}^{-1}     \\
        &=&
        \tilde{k}\tilde{h}\left( n\cartan_{(k,y)}^\ast n^{-1} \right) \tilde{h}^{-1}\tilde{k}^{-1}
                                                                    \\
        &=&
        \tilde{k}\tilde{h}n
        \left( \tilde{h}^{-1}\tilde{k}^{-1}\cartan_{(k^\prime, y^\prime)}^\ast\tilde{k}
            \tilde{h} \right) n^{-1} \tilde{h}^{-1}\tilde{k}^{-1}.
\end{array}
\]
Using the fact that $n \in N_H(\cartan_{(k,y)})$, a routing computation verifies that
$m:= \tilde{k}\tilde{h}n \tilde{h}^{-1}\tilde{k}^{-1} \in N_H(\cartan_{(k^\prime, y^\prime)})$.
But then $\tilde{h}(l, z) \in (V_{(h,x)})_{K^\prime}$,
and $\cartan_{\tilde{h}(l, z)}^\ast = m\cartan_{(k^\prime, y^\prime)}^\ast m^{-1}$
with $m \in N_H(\cartan_{(k^\prime, y^\prime)})$.  It follows that
$\tilde{h}(l, z) \in \mathcal{U}_{\tilde{g}}^{\cartan_{(k^\prime, y^\prime)}}(k^\prime, y^\prime)$.
Since $\mathcal{U}_{\tilde{g}}^{\cartan_{(k,y)}}(\tilde{k}, \tilde{y})$ is $G$-invariant,
\[
    \mathcal{U}_{\tilde{g}}^{\cartan_{(k,y)}}(\tilde{k}, \tilde{y})
    \subseteq
    \mathcal{U}_g^{\cartan_{(k^\prime, y^\prime)}}(k^\prime, y^\prime).
\]
Switching the roles of $(k, y)$ and $(k^\prime, y^\prime)$ completes the proof.
\end{proof}

Note that we may now denote $\mathcal{U}_{\tilde{g}}^{\cartan_{(k,y)}}(\tilde{k}, \tilde{y})$
simply as $\mathcal{U}^{\cartan_{(k,y)}}(\tilde{k}, \tilde{y})$.

\begin{lemma}
\label{lem-StratConnectedPieceIndepPoint}
If $(\tilde{l}, \tilde{z}) \in \mathcal{U}^{\cartan_{(k,y)}}(k,y)$, $(l,z) \in V_{(h,x)}$ is
in the same orbit as $(\tilde{l}, \tilde{z})$, and $\cartan_{(l, z)}^\prime$  a Cartan subgroup
of $G_{(l, z)}$ associated to $l$, then
$\mathcal{U}^{\cartan_{(k,y)}}( \tilde{k},\tilde{y}) = \mathcal{U}^{\cartan_{(l, z)}^\prime}(\tilde{l},\tilde{z})$.
In particular, $\mathcal{U}^{\cartan_{(k,y)}}(\tilde{k},\tilde{y})$ does not depend on the choice
of a Cartan subgroup $\cartan_{(k,y)}$ of $K$ associated to $k$.
\end{lemma}
\begin{proof}
Let $K = Z_{G_y}(k)$ as above. By Eq.~\eqref{eq:decUrep}, we may assume that
$(l, z) \in (V_{(h,x)})_K \cap (\cartan_{(k,y)}\times M)$, and that
$\cartan_{(l, z)}^\ast = n \cartan_{(k,y)}^\ast n^{-1}$ for some $n \in N_H (\cartan_{(k,y)})$.
Let $\cartan_{(l, z)}^\prime$ be a choice of a Cartan subgroup
in $K$ associated to $l$. Then, by \cite[IV.~Prop.~4.6]{BroeckertomDieck}, there is an
$i \in K$ such that $i\cartan_{(l, z)}^\prime i^{-1} = \cartan_{(k,y)}$.
Recall that by Lemma \ref{lem-StratConnectedTorusinG},
$i\cartan_{(l, z)}^{\prime\,\ast} i^{-1} = \cartan_{(l, z)}^\ast$.

Given $(j, w) \in (V_{(h,x)})_K \cap (\cartan_{(l, z)}^\prime \times M)$ such that
$\cartan_{(j, w)}^{\prime\,\ast} = m\cartan_{(l, z)}^{\prime\,\ast} m^{-1}$
for some $m \in N_H(\cartan_{(l, z)}^\prime)$, it is now easy to see that
\[
    (j, w) =
    i(j, w)
    \in
    i \left( (V_{(h,x)})_K \cap (\cartan_{(l, z)}^\prime \times M)\right)
    =
    (V_{(h,x)})_K \cap (\cartan_{(k,y)} \times M).
\]
In particular, $j \in \cartan_{(k,y)}$ so that by Lemma \ref{lem-StratConnectedTorusinG},
$\cartan_{(j, w)}^\ast = i\cartan_{(j, w)}^{\prime\,\ast} i^{-1}$.
Then
\[
    \cartan_{(j, w)}^\ast
    =   i\cartan_{(j, w)}^{\prime\,\ast} i^{-1}
    =   imi^{-1}n\cartan_{(k,y)}^\ast n^{-1}i m^{-1}i^{-1}.
\]
A routine computation
verifies that $imi^{-1}n \in N_H(\cartan_{(k,y)})$.
Therefore, $(j, w) \in \mathcal{U}^{\cartan_{(k,y)}}(\tilde{k},\tilde{y})$, so that
$\mathcal{U}^{\cartan_{(l, z)}^\prime}(\tilde{l},\tilde{z}) \subseteq \mathcal{U}^{\cartan_{(k,y)}}(\tilde{k},\tilde{y})$.
Switching the roles of $(\tilde{k},\tilde{y})$ and $(\tilde{l},\tilde{z})$ completes the proof that
$\mathcal{U}^{\cartan_{(l, z)}^\prime}(\tilde{l},\tilde{z}) = \mathcal{U}^{\cartan_{(k,y)}}(\tilde{k},\tilde{y})$.

If $\cartan_{(k,y)}^\prime$ is another choice of a Cartan subgroup of $K$
associated to $k$,
repeating the above argument with $(l, z) = (k,y)$ yields
$\mathcal{U}^{\cartan_{(k,y)}}(k,y) = \mathcal{U}^{\cartan_{(k,y)}^\prime}(k,y)$.
\end{proof}

Because of the preceding considerations, the set $\mathcal{U}^{\cartan_{(k,y)}}(\tilde{k},\tilde{y})$
depends only on the orbit $G(k,y)$, hence we will denote it simply as $\mathcal{U}\big( G(k,y) \big)$.
For $(\tilde{k} , \tilde{y})$ in the same orbit as $(k,y)$, we denote by
$\mathcal{U}\big( G(k,y) \big)^\textrm{c}_{(\tilde{k} , \tilde{y})}$ or even shorter by
$\mathcal{U}_{(\tilde{k} , \tilde{y})}$ the connected component of
$(\tilde{k} , \tilde{y})$ in $\mathcal{U}\big( G(k,y) \big)$. The partition $\mathcal{Z}$
of $U$ then can be written as
\begin{equation}
\label{eq-StratConnectedPieces}
 \mathcal{Z} = \big\{ \mathcal{U}_{(\tilde{k} , \tilde{y})} \in \mathcal{P} (U)\mid
             (\tilde{k} , \tilde{y}) \in U \big\} .
\end{equation}

Having established that the sets $\mathcal{U}\big( G(k,y) \big)$ are well-defined,
we now confirm that the set germs of the $\mathcal{U}\big( G(k,y) \big)$ coincide with
the stratification given by Equation \eqref{eq-StratTildeMDef}.

\begin{proposition}
\label{prop-StratConnectedGermsCoincide}
For each $(\tilde{k}, \tilde{y}) \in U$,
the germs $[\mathcal{U}\big( G(k,y) \big)]_{(\tilde{k},\tilde{y})}$,
$[\mathcal{U}_{(\tilde{k},\tilde{y})}]_{(\tilde{k},\tilde{y})}$ and $\mathcal{S}_{(\tilde{k},\tilde{y})}$ coincide.
\end{proposition}
\begin{proof}
As $\mathcal{S}_{(\tilde{k}, \tilde{y})}$ and $\mathcal{U}\big(G (k,y)\big)$
depend only on the orbit of $(\tilde{k}, \tilde{y})$, it is clearly
sufficient to consider the case of $(\tilde{k}, \tilde{y}) = (k, y) \in V_{(h,x)}$.
Set $K = Z_{G_y}(k)$ and fix a linear slice $V_{(k,y)}$ at $(k, y)$ for
the $G$-action on $G \times M$.  By \cite[II. Corollary 4.6]{BredonBook},
we may assume that $V_{(k,y)} \subseteq V_{(h,x)}$,
though it need not be the case that $V_{(k,y)}$ is the image under the exponential map
of a subset of the normal space $V_{(k,y)}$.  As in the proof of Lemma
\ref{lem-StratConnectedSliceIndependent}, we define a closed subset $C$ of $G$ consisting of
the (finitely many) connected components of $K$ not containing $k$ as well as the (finitely many)
nontrivial $N_H(\cartan_{(k,y)})$-conjugates of $\cartan_{(k,y)}^\ast$.  Let $O = C^c$
be the complement of $C$ in $G$. Then $O \times M$ is an open
subset of $G \times M$ containing $\cartan_{(k,y)}^\ast$.
Hence $V_{(k,y)} \cap (O\times M)$ is an open neighborhood of $(k, y)$
in $V_{(k,y)}$, so we may shrink $V_{(k,y)}$ to assume that $V_{(k,y)} \subseteq O\times M$.
Put $Q = GV_{(k,y)}$. We now show that the set germs $[\mathcal{U}\big( G(k,y) \big)]_{(k,y)}$ and $\mathcal{S}_{(k,y)}$
coincide by proving that
\[
    \mathcal{U}\big( G(k,y) \big) \cap Q
    =
    G\left( V_{(k,y)}^K \cap (\cartan_{(k,y)}^\ast\times M)\right).
\]

Let $(\tilde{l}, \tilde{z}) \in \mathcal{U}\big( G(k,y) \big) \cap Q$.  Then there is a
$\tilde{g}^\prime \in G$ such that
$(l, z) := \tilde{g}^\prime(\tilde{l}, \tilde{z})
\in (V_{(h,x)})_K \cap (\cartan_{(k,y)} \times M)$
and an $n \in N_H(\cartan_{(k,y)})$ such that
$\cartan_{(l, z)}^\ast = n\cartan_{(k,y)}^\ast n^{-1}$.
In particular, $\cartan_{(k,y)}$ is a Cartan subgroup of $K$ associated to $l$.
As $(l, z) \in Q$, there is a $g \in G$ such that $g(l, z) \in V_{(k,y)}$.
Moreover, as $G_{g(l, z)} \leq K$ and $G_{(l, z)} \leq K$, we have $G_{g(l, z)} = K$
by \cite[Lem.~4.2.9]{PflaumBook} and hence $g \in N_H(K)$.
Similarly, as $g(l, z) \in V_{(k,y)} \subseteq V_{(h,x)}$ and $(l, z) \in V_{(h,x)}$,
$g \in H$ by (SL4).

As $k, glg^{-1} \in K$ and $(k,y), g(l, z) \in V_{(k,y)}$, which is disjoint from
$(K\smallsetminus kK^\circ)\times M$, $k$ and $glg^{-1}$ are in the same connected
component of $K$.  By \cite[IV.~Prop.~4.6]{BroeckertomDieck}, there is a
$\tilde{k} \in K$ such that $\tilde{k}g\cartan_{(k,y)}g^{-1}\tilde{k}^{-1} = \cartan_{(k,y)}$,
and hence $\tilde{k}g \in N_H(\cartan_{(k,y)})$.
Recalling that $l \in \cartan_{(l, z)}^\ast = n\cartan_{(k,y)}^\ast n^{-1}$
for some $n \in N_H(\cartan_{(k,y)})$, we have
\[
    \tilde{k}g l g^{-1} \tilde{k}^{-1}
    \in
    \tilde{k}g\cartan_{(l, z)}^\ast g^{-1}\tilde{k}^{-1}
    =
    \tilde{k}gn\cartan_{(k,y)}^\ast n^{-1}g^{-1}\tilde{k}^{-1}.
\]
Recalling that $\tilde{k} \in K$, and $K$ is the isotropy group of $g(l, z) = (gz, glg^{-1})$, we have
\[
    g l g^{-1}
    \in
    \tilde{k}gn\cartan_{(k,y)}^\ast n^{-1}g^{-1}\tilde{k}^{-1}.
\]
However, as $g(l, z) \in V_{(k,y)}\subseteq O \times M$,
which is disjoint from $C \times M$, and as $\tilde{k}gn \in N_H(\cartan_{(k,y)})$, it must be that
$\tilde{k}gn\cartan_{(k,y)}^\ast n^{-1}g^{-1}\tilde{k}^{-1} = \cartan_{(k,y)}^\ast$.
It follows that
$g(l, z) \in V_{(k,y)}^K \cap (\cartan_{(k,y)}^\ast\times M)$,
and hence $\mathcal{U}\big( G(k,y) \big) \cap Q \subseteq
G\left( V_{(k,y)}^K \cap (\cartan_{(k,y)}^\ast\times M)\right)$.

Conversely, if
$(\tilde{j}, \tilde{w}) \in G\left( V_{(k,y)}^K \cap (\cartan_{(k,y)}^\ast\times M)\right)$,
then there is a $\hat{g} \in G$ such that
$(j, w) := \hat{g}(\tilde{j}, \tilde{w}) \in V_{(k,y)}^K \cap (\cartan_{(k,y)}^\ast\times M)$.
Then as $V_{(k,y)} \subseteq V_{(h,x)}$, we have
\[
    (j, w) \in (V_{(h,x)})_K \cap (\cartan_{(k,y)}^\ast\times M),
\]
and so $(j, w) \in \mathcal{U}\big( G(k,y) \big)$ using the trivial element of the normalizer.
Therefore,
\[
    \mathcal{U}\big( G(k,y) \big) \cap Q
    =
    G\left( V_{(k,y)}^K \cap (\cartan_{(k,y)}^\ast\times M)\right),
\]
which shows the first part of the claim.

As the set $\left( V_{(k,y)}^K \cap (\cartan_{(k,y)}^\ast\times M)\right) $ in the right hand side of the
preceding equation is connected, the set $\mathcal{U}\big( G(k,y) \big) \cap Q$ has a finite number of
connected components. Since $Q$ is a $G$-invariant open neighborhood of the orbit $G(k,y)$, this
implies that $\mathcal{U}\big( G(k,y) \big)$ has locally only finitely many connected components
and that the germ of $\mathcal{U}_{(k,y)}$ at $(k,y)$ coincides
with the germ of $\mathcal{S}_{(k,y)}$ at $(k,y)$.
\end{proof}

Since the $\mathcal{S}_{(h,x)}$ are germs of smooth $G$-submanifolds of
$G \times M$, and the piece associated to a point $(\tilde{k},\tilde{y}) \in U$ has the same
set germ as $\mathcal{S}_{(l, z)}$ at $(l, z) \in \mathcal{U}_{(\tilde{k},\tilde{y})}$,
it follows that the pieces of $\mathcal{Z}$ are smooth submanifolds of $G \times M$
invariant under the $G$-action.

\begin{lemma}
\label{lem-StratConnectedLocFinite}
The partition $\mathcal{Z}$ of $U = GV_{(h,x)}$ given by Equation \eqref{eq-StratConnectedPieces}
is finite.
\end{lemma}
\begin{proof}
The set $\mathcal{U}\big( G(k,y) \big)$ is
determined by the $H$-conjugacy class of $G_{(k,y)} = Z_{G_y}(k) \leq H$
as well as the connected component containing $k$ of the $\sim$ class of $k$ for the
$\cartan_{(k,y)}$-action on $GV_{(k,y)}$.
By Lemma \ref{lem-StratConnectedIndepTorus}, the set $\mathcal{U}\big( G(k,y) \big)$ does not
depend on the choice of a Cartan subgroup associated to $k$.
As $H$ acts linearly on $N_{(h,x)}$,
there is a finite number of $H$-conjugacy classes of isotropy
groups with respect to the $H$-action on  $N_{(h,x)}$, hence on $V_{(h,x)}$.

Choose a representative $K$ of each $H$-isotropy type in $V_{(h,x)}$.  Then as
$K/K^\circ$ is finite, there are a finite number of conjugacy classes of cyclic
subgroups of $K/K^\circ$ and hence by \cite[IV.~Prop.~4.6]{BroeckertomDieck}
a finite number of $K$-conjugacy classes of Cartan subgroups of $K$.  Given
a Cartan subgroup $\cartan_{(k,y)}$ of $K$, there are a finite number of connected components
of $\sim$ classes in $\cartan_{(k,y)}$ with respect to the action of $\cartan_{(k,y)}$ on
$U$.  Of course, $\cartan_{(k,y)}^\ast$ is defined with respect to the action of
$\cartan_{(k,y)}$ on a subset of $U$, but this implies that $\cartan_{(k,y)}^\ast$
is given by a union of connected components of $\sim$ classes with respect to the action on $U$,
of which there are finitely many.  It follows that there are a finite number of
$\mathcal{U}\big( G(k,y) \big)$.  Finally, each $\mathcal{U}\big( G(k,y) \big)$ has a finite number of connected components,
which completes the proof.
\end{proof}

We now verify that $\mathcal{Z}$ is a decomposition indeed, cf.~\cite[Def.~1.1.1 (DS2)]{PflaumBook}.

\begin{proposition}
\label{prop-StratConnectedFrontier}
The pieces of $\mathcal{Z}$ satisfy the condition of frontier.
\end{proposition}
\begin{proof}
Suppose $\mathcal{U}\big( G(k,y) \big) \cap \overline{\mathcal{U}\big( G(l,z) \big)} \neq \emptyset$
where the closure is taken in $\inertia{M}$.  As the pieces of $\mathcal{Z}$ are defined to be the
connected components of the $\mathcal{U}\big( G(k,y) \big)$ and $\mathcal{U}\big( G(l,z) \big)$,
it is sufficient to show that
$\mathcal{U}\big( G(k,y) \big) \cap \overline{\mathcal{U}\big( G(l,z) \big)}$ is both open
and closed in $\mathcal{U}\big( G(k,y) \big)$.  It is obvious that
$\mathcal{U}\big( G(k,y) \big) \cap \overline{\mathcal{U}\big( G(l,z) \big)}$ is closed in
$\mathcal{U}\big( G(k,y) \big)$, so we need only establish that
$\mathcal{U}\big( G(k,y) \big) \cap \overline{\mathcal{U}\big( G(l,z) \big)}$ is open in
$\mathcal{U}\big( G(k,y) \big)$.

By Lemma \ref{lem-StratConnectedPieceIndepPoint}, the piece $\mathcal{U}\big( G(k,y) \big)$
may be defined in terms of any orbit it contains, so we may assume that some element of the $G$-orbit
of $(k,y)$ is contained in $\mathcal{U}\big( G(k,y) \big) \cap \overline{\mathcal{U}\big( G(l,z) \big)}$.
Then the $G$-invariance of these two sets implies that
$G(k,y) \subseteq \mathcal{U}\big( G(k,y) \big) \cap \overline{\mathcal{U}\big( G(l,z) \big)}$.
By Proposition \ref{prop-StratConnectedGermsCoincide}, an open neighborhood of
$(k, y)$ in $\mathcal{U}\big( G(k,y) \big)$ is given by
$G\big(V_{(k, y)}^K \cap (\cartan_{(k, y)}^\ast\times M)\big)$
for a sufficiently small slice $V_{(k, y)}$ at $(k, y)$.  As above, we may assume
$V_{(k, y)} \subseteq V_{(h,x)}$ by \cite[II. Corollary 4.6]{BredonBook} so that
while $V_{(k, y)}$ can then be taken to be linear, it need not be the image under the
exponential map of a subset of the normal space $N_{(k,y)}$.
We will show that $G\big(V_{(k, y)}^K \cap (\cartan_{(k, y)}^\ast\times M)\big)$
is contained in $\overline{\mathcal{U}\big( G(l,z) \big)}$.

As $GV_{(k, y)}$ must contain some element of $\mathcal{U}\big( G(l,z) \big)$, we may assume again
by Lemma \ref{lem-StratConnectedPieceIndepPoint} that $G(l,z) \subseteq GV_{(k,y)}$.  Moreover,
by the proof of Lemma \ref{lem-StratConnectedUkyFiniteCC}, we may choose the representative
$(l,z)$ from the orbit $G(l,z)$ such that $(l,z) \in V_{(k,y)}$,
$k \in \cartan_{(l,z)} \leq \cartan_{(k,y)}$, and $k \in \overline{\cartan_{(l,z)}^\ast}$.
Let $K = G_{(k,y)}$ and $L = G_{(l,z)}$ so that $L \leq K$,
and then $(k,y) \in V_{(k,y)}^K \subseteq (V_{(h,x)})_K \subseteq \overline{(V_{(h,x)})_L}$.
Then we have
\[
    (k,y) \in \overline{(V_{(h,x)})_L}
    \cap \big(\overline{\cartan_{(l,z)}^\ast}  \times M \big).
\]
In particular, note that by our choice of $(l,z) \in V_{(k,y)}$ used to define the set
$\mathcal{U}\big( G(l,z) \big)$, $(k,y)$ is in the closure of the set corresponding to
the trivial element of $N_H(\cartan_{(l,z)})$ in Equation \eqref{eq:decUrep}.

For any $(j, w) \in V_{(k, y)}^K \cap (\cartan_{(k, y)}^\ast\times M)$,
as $j \in \cartan_{(k, y)}^\ast$, it follows that $(GV_{(k,y)})^j = (GV_{(k,y)})^k$.
In particular, $k \in \cartan_{(l,z)} \leq L$ implies that $k$ fixes $(l, z) \in V_{(k,y)}$ so that
$j \in L$ as well.  Since $V_{(k, y)}^K \cap (\cartan_{(k, y)}^\ast\times M)$ is invariant under the action of scalar
$\in [0,1]$,
and $k$ is in the same connected component of $L$ as $l$,  each such $j$ is in
the same connected component of $L$ as $l$ also.
Fix a $(j,w) \in V_{(k, y)}^K \cap (\cartan_{(k, y)}^\ast\times M)$. Then
there is a $\tilde{l} \in L$ such that
$\tilde{l}\cartan_{(l,z)}\tilde{l}^{-1}$ is a Cartan subgroup of $L$ associated to $j$.
Hence $\tilde{l}^{-1} j \tilde{l} \in \cartan_{(l,z)}$, so that as $\tilde{l} \in L \leq K = G_{(j,w)}$,
we have $\tilde{l}^{-1} j \tilde{l} = j \in \cartan_{(l,z)}$.

Finally, note that as $j \in \cartan_{(k, y)}^\ast$, it is clear that
$(GV_{(l,z)})^j = (GV_{(l,z)})^k$ for a slice $V_{(l,z)}$ chosen small enough so that
$GV_{(l,z)} \subseteq GV_{(k,y)}$.  Therefore $j \sim k$ as elements of $\cartan_{(l,z)}$.
Then as the connected component $[k]^\circ$ of the $\sim$ class of $k$ as an
element of $\cartan_{(l,z)}$ intersects $\overline{\cartan_{(l,z)}^\ast}$,
we have by Lemma \ref{lem-StratConnectedTorusFrontier} that
$[k]^\circ \subseteq \overline{\cartan_{(l,z)}^\ast}$.
It follows that
\[
    (j,w) \in \overline{(V_{(h,x)})_L}
    \cap \big( \overline{\cartan_{(l,z)}^\ast}  \times M \big),
\]
so as $(j,w) \in V_{(k, y)}^K \cap (\cartan_{(k, y)}^\ast\times M)$ was arbitrary,
\[
    V_{(k, y)}^K \cap (\cartan_{(k, y)}^\ast\times M)
    \subseteq \overline{(V_{(h,x)})_L}
    \cap \big(\overline{\cartan_{(l,z)}^\ast} \times M \big).
\]
Considering the $G$-saturations of both sides of this inclusion, it follows that each element of
$\mathcal{U}\big( G(k,y) \big) \cap \overline{\mathcal{U}\big( G(l,z) \big)}$ is contained
in a neighborhood that is both open and closed in $\mathcal{U}\big( G(k,y) \big)$,
completing the proof.
\end{proof}

Finally, we have the following.

\begin{proposition}
\label{prop-StratConnectedWhitneyB}
The orbit Cartan type stratifications of $\inertia{M}$ and the inertia space
$\inertia{X}$ both satisfy Whitney's condition B.
\end{proposition}

The proof follows \cite[Thm.~4.3.7]{PflaumBook}.

\begin{proof}
Let $(h,x) \in \inertia{M}$,  $H = Z_{G_x}(h)$, and $V_{(h,x)}$ a slice
at $(h,x)$ of the form $\exp (B_{(h,x)})$, where $B_{(h,x)}$ is a ball around the origin
in the normal space $N_{(h,x)}$. We work in the neighborhood $U:=GV_{(h,x)}$
of $(h,x)$ in $\inertia{M}$, and show that for any stratum $S \in \mathcal{Z}$
with $(h,x)\in \overline{S}$ Whitney's condition B is satisfied  at $(h,x)$ for
the pair of strata $(R,S)$,
where $R$ is the piece of $\mathcal{Z}$ containing $(h,x)$.
Recall that $\mathcal{Z}$ is the decomposition of $U$ given by
Eq.~\eqref{eq-StratConnectedPieces}. Recall also, that $R$ is the
connected component of $G\left( V_{(h,x)}^H \cap (\cartan_{(h,x)}^\ast\times M)\right)$
containing $(h,x)$. To describe the stratum $S$ in some more detail, consider
an orbit $G(k,y)$ for $(k, y) \in S$. As in the proof of Lemma \ref{lem-StratConnectedUkyFiniteCC},
we may choose the representative $(k,y)$ of the orbit $G(k,y)$ such that
$(k,y) \in V_{(h,x)}$, $h \in \cartan_{(k,y)} \leq \cartan_{(h,x)}$, and
$h \in \overline{\cartan_{(k,y)}^\ast}$. In particular, we then have the relation
$K  \leq H$ for the isotropy group $K := Z_{G_y}(k)$ of $(k,y)$.
As shown above, $S$  coincides with the
connected component of $\mathcal{U}\big( G(k,y) \big)$ containing
$(k,y)$.

Suppose now that
$(h_i, x_i)_{i \in \N}$ is a sequence in $R$ and $(k_i, y_i)_{i \in \N}$  a sequence in $S$,
and that both sequences converge to $(h,x)$.
Assume in addition that in a smooth chart around $(h,x)$ the secants
$\ell_i = \overline{(h_i,x_i),(k_i,y_i)}$ converge to a straight line $\ell$, and the
tangent spaces $T_{(k_i,y_i)}S$ converge to a subspace $\tau$.
Then we must show that $\ell \subseteq \tau$.

Note that the hypotheses imply that
$(h,x) \in \mathcal{U}\big( G(h,x) \big)\cap \overline{\mathcal{U}\big( G(k,y) \big)}$.
By the proof of Proposition~\ref{prop-StratConnectedFrontier}
and the choices of $(k,y)$ and $\cartan_{(k,y)} \subseteq K$
we obtain the relation
\begin{equation}
\label{eq:closurerel1} V_{(h, x)}^H \cap (\cartan_{(h, x)}^\ast\times M)
    \subseteq \overline{(V_{(h,x)})_K}
    \cap \big(\overline{\cartan_{(k,y)}^\ast} \times M \big) .
\end{equation}
Moreover, since  every element $n \in N_H(\cartan_{(k,y)})$ fixes
$V_{(h, x)}^H \cap (\cartan_{(h, x)}^\ast\times M)$, it follows
that $V_{(h, x)}^H \cap (\cartan_{(h, x)}^\ast\times M)
\subseteq n\overline{\cartan_{(k,y)}^\ast}n^{-1} \times M$ as well, hence
\begin{equation}
\label{eq:closurerel2} V_{(h, x)}^H \cap (\cartan_{(h, x)}^\ast\times M)
    \subseteq \overline{(V_{(h,x)})_K}
    \cap \big(n \overline{\cartan_{(k,y)}^\ast} n^{-1} \times M \big) .
\end{equation}

Denote by $\mathfrak{g}$ the Lie algebra of $G$, by $\mathfrak{h}$ the Lie algebra
of $H$, and let  $\mathfrak{m}$ denote the orthogonal complement of $\mathfrak{h}$ in
$\mathfrak{g}$ with respect to the initially chosen bi-invariant metric on $G$.
Then there is a neighborhood $U^\prime \subseteq U \cong G \times_H V_{(h,x)}$ of $(h,x)$
in $G \times M$ such that
\[
    \Psi :   U^\prime   \longrightarrow  \mathfrak{m} \times N_{(h,x)} , \:
    [\exp_{|\mathfrak{m}} \xi, \exp_{(h,x)}(v)]\longmapsto (\xi, v)
\]
is a smooth chart at $(h,x)$, where $\exp_{|\mathfrak{m}}$ denotes the restriction of the exponential map of the Lie group $G$ to $\mathfrak{m}$, and
$\exp_{(h,x)}$ the exponential function restricted to the open ball
$B_{(h,x)} \subseteq N_{(h,x)}$.
By shrinking $U^\prime$ if necessary, we have that there is an open neighborhood $O$
of $H$ in $G$ such that
\[
    \Psi\left(O\left( V_{(h,x)}^H \cap (\cartan_{(h,x)}^\ast\times M)\right) \right)
    \subset
    \mathfrak{m} \times \left( N_{(h,x)}^H \cap T_{(h,x)} (\cartan_{(h,x)}^\ast\times M)\right).
\]
We may assume that the sequences
$(h_i, x_i)_{i \in \N}$ and $(k_i, y_i)_{i \in \N}$ are contained in $U^\prime$.
Since $(k_i, y_i) \in \mathcal{U}\big( G(k,y) \big)$, one knows that
\[
    \Psi(k_i,y_i)
    \in \mathfrak{m} \times H \left(
    (V_{(h,x)})_K \cap \left(\bigcup\limits_{n \in N_H (\cartan_{(k,y)})}
                n\cartan_{(k,y)}^\ast n^{-1}\times M\right)\right).
\]
Recall that there are only finitely many and pairwise disjoint
sets $n\cartan_{(k,y)}^\ast n^{-1}$, where $n$ runs through the elements of
$N_H (\cartan_{(k,y)})$. Moreover, by Lemma
\ref{lem-StratConnectedNormalizerPermut}, $n\cartan_{(k,y)}^\ast n^{-1}$  is disjoint from
$\overline{m\cartan_{(k,y)}^\ast m^{-1}} = m \left(\overline{\cartan_{(k,y)}^\ast}\right) m^{-1}$ for every $m \in N_H (\cartan_{(k,y)})$ with
$n\cartan_{(k,y)}^\ast n^{-1} \neq m\cartan_{(k,y)}^\ast m^{-1}$.
Hence, we may assume without loss of generality that
\[
    (k_i, y_i) \in G \left((V_{(h,x)})_K \cap
    \left( m_0 \cartan_{(k,y)}^\ast m_0^{-1}\times M\right)\right)
\]
for all $i$ and some $m_0 \in  N_H (\cartan_{(k,y)})$.

Choose $\tilde{l}_i \in G$ such that
$(\tilde{k}_i, \tilde{y}_i):= \tilde{l}_i (k_i, y_i) \in (V_{(h,x)})_K$
for all $i\in \N$. Put $(\tilde{h}_i,\tilde{x}_i):= l_i (h_i,x_i)$. After possibly
passing to a subsequence, $(\tilde{l}_i)_{i\in \N}$
converges to some $\tilde{l} \in H$,  the secant lines $\tilde{\ell}_i =
\overline{(\tilde{h}_i,\tilde{x}_i),(\tilde{k}_i,\tilde{y}_i)}$ converge to a straight
line $\tilde{\ell}$, and the
tangent spaces $T_{(\tilde{k}_i,\tilde{y}_i)}S$ converge to a subspace $\tilde{\tau}$.
By definition, and since
$\tilde{l}_i T_{(k_i,y_i)} S = T_{(\tilde{k}_i,\tilde{y}_i)}S$ for all $i$,
one obtains $\tilde{\ell} = \tilde{l}  \ell$, and $\tilde{\tau} = \tilde{l} \tau$.
Hence, the first claim is shown, if
$\tilde{\ell} \subseteq \tilde{\tau}$. Without loss of generality we may therefore
assume that for all $i \in \N$
\begin{equation}
\label{eq:reductionky}
  (k_i, y_i) \in  (V_{(h,x)})_K \cap \big( m_0\cartan_{(k,y)}^\ast m_0^{-1}
        \times M \big) ,
\end{equation}
and then show $\ell \subseteq \tau$ for the sequences $(k_i, y_i)_{i\in \N}$
and $(h_i,x_i)_{i\in  \N}$.

Eq.~\eqref{eq:reductionky} now means in particular that
\[
    \Psi (k_i,y_i)
    \in \{ 0 \} \times \Big(
    (N_{(h,x)})_{K} \cap
    \exp^{-1}_{(h,x)}\big( m_0\cartan_{(k,y)}^\ast m_0^{-1}\times M\big)\Big).
\]
Since by Lemma \ref{lem-StratConnectedTorusPartition}
and the above observations $m_0 \overline{\cartan_{(k,y)}^\ast} m_0^{-1}$ is an open
and closed subset of a closed subgroup of $G$ and also contains $h$, the set
\[
    V:= N_{(h,x)} \cap T_{(h,x)}
    \Big( m_0 \big(\overline{\cartan_{(k,y)}^\ast}\big) m_0^{-1}\times M\Big)
\]
is a subspace of $N_{(h,x)}$. Let $W$ be the orthogonal complement of the invariant
space $V^H$ in $V$ with respect to the $H$-invariant scalar product induced from
$V_{(h,x)}$. Then the image  under the chart $\Psi$ of every element of
$G\big( V_{(h,x)}^H \cap (\cartan_{(h,x)}^\ast\times M)\big)\cap U^\prime$
and every $(k_i, y_i)$ is contained in
\[
    \mathfrak{m} \times (W_K \cup \{ 0 \}) \times V^H .
\]
With respect to this decomposition, $(h,x)$ has coordinates $(0, 0, 0)$,
each element of $G\left( V_{(h,x)}^H \cap (\cartan_{(h,x)}^\ast\times M)\right)$ has
coordinates contained in $\mathfrak{m} \times 0 \times V^H$, and each sequence
element $(k_i, y_i)$ has coordinates contained in
$\{ 0 \}  \times W_{K} \times V^H$.  In particular,
let
\[
    \Psi(k_i,y_i)
    =
    (0 , w_i, v_i)
\]
for every $i$.  Then as $W_{K}$ is invariant under multiplication by
non-vanishing scalars, we have
\[
\begin{split}
    (\xi, w, v)
     :=\, &
    \lim\limits_{i\to\infty} \frac{ \Psi(k_i,y_i) - \Psi(h_i,x_i) }
        {\| \Psi(k_i,y_i) - \Psi(h_i,x_i) \|}
    \, \in \mathfrak{m} \times \overline{W_{K}} \times V^H
\end{split}
\]
Now, as the unit sphere in $W$ is compact, the sequence $\frac{w_i}{\| w_i \|}$
converges to some $\hat{w} \in SW$ after possibly passing to a subsequence.
Then $w =  \| w \|\hat{w}$.
Since $W_{K}$ is invariant by non-vanishing scalars, we have
\[
    \mathfrak{m} \times \operatorname{span}\: \hat{w} \times V^H \subseteq \tau,
\]
and
\[
    \ell = \operatorname{span} \: (\xi, \hat{w}, v) \subseteq \tau,
\]
proving the first claim.

Now let us show that the orbit Cartan type stratification of $\inertia{X}$
satisfies Whitney's condition B as well.
To this end let us first choose a Hilbert basis of $H$-invariant polynomials
$p_1,\ldots , p_\kappa : \big(N_{(h,x)}^H\big)^\perp \rightarrow \R$ of the orthogonal complement of
the invariant space $N_{(h,x)}^H$ in $N_{(h,x)}$.
Next let $p_{\kappa + 1},\ldots , p_N : N_{(h,x)}^H \rightarrow \R$ with $N=\kappa +\dim N_{(h,x)}^H$ be a linear
coordinate system of the invariant space. We can even choose these $p_i$ in such a way
that $p_{\kappa + 1},\ldots , p_{\kappa + \dim V^H}$ is a linear coordinate
system of $V^H$. By construction, $p_1,\ldots , p_N$ then is a Hilbert basis of the normal space $N_{(h,x)}$.
Denote by $p: N_{(h,x)} \rightarrow \R^N$ the corresponding Hilbert map.
Recall that $p$ induces a chart of $\inertia{X} $ over
$ G \backslash U$ by
\[
  \widehat{\Psi}: G \backslash U \rightarrow \R^N, \: G \exp_{(h,x)} (v) \mapsto p(v) .
\]
Note that by $H$-invariance of $p$ and since for every orbit in $U$ there is
a representative in $V_{(h,x)}$, the chart $\widehat{\Psi}$ is well-defined indeed.
A decomposition of $\widehat{U}:= \widehat{\Psi} ( G \backslash U )$ inducing the
orbit Cartan type stratification on $G \backslash U$  is given by
\[
  \widehat{\mathcal{Z}}:= \big\{ \widehat{\Psi} (G\backslash GS) \mid S \in \mathcal{Z} \big\} .
\]
Let $\widehat{S}\in \widehat{\mathcal{Z}}$ denote the stratum containing the orbit
$G(h,x)$, and $\widehat{S}\in \mathcal{Z}$ a stratum $\neq \widehat{R}$ such that
$G(h,x)$ lies in the closure of $\widehat{S}$.
Now consider sequences of orbits $\big( G(h_i,x_i) \big)_{i\in \N}$ in $\widehat{R}$
and $\big( G(k_i,y_i) \big)_{i\in \N}$ in  $\widehat{S}$ such that both sequences converge
to $G(h,x)$. Moreover, assume  that the sequence of secants
$ \overline{\widehat{\Psi} (G(h_i,x_i)) ,  \widehat{\Psi} (G(k_i,y_i))}$
converges to a line $\widehat{\ell}$, and that the sequence of tangent spaces
$T_{\widehat{\Psi} (G(k_i,y_i))} \widehat{S}$ converges to some subspace
$\widehat{\tau} \subseteq \R^N$. Using notation from before, we can
choose representatives $(h_i,x_i)$ and $(k_i,y_i)$  having coordinates in
$ \mathfrak{m} \times (W_K \cup \{ 0 \}) \times V^H \subseteq N_{(h,x)}$ such that
\begin{equation}
\label{eq:coordinates}
\begin{split}
  \Psi (h_i,x_i) & \, = (0,0, v_i^\prime ) \in \{ 0 \}  \times  \{ 0 \} \times V^H
  \: \text{ and } \\
  \Psi (k_i,y_i) & \, = (0,w_i,v_i) \in \{ 0 \} \times W_K  \times V^H .
\end{split}
\end{equation}

Next observe that by the Tarski--Seidenberg Theorem and the proof of Lemma
\ref{lem-StratConnectedUkyFiniteCC}, the stratum $\widehat{S}$ is semialgebraic
as the image of the semialgebraic set $(W_K \times V^H) \cap B_{(h,x)}$ under the Hilbert map $p$.
By the same argument, $p(W_K)$ is semialgebraic, too, and an analytic manifold, since
$p(W_K) \cong N_H(K)\backslash W_K \cong H \backslash W_{(K)}$. Moreover, the equality
\[
  \widehat{S} = ( p(W_K) \times V^H) \cap p(B_{(h,x)} )
\]
holds true, where we have canonically identified $V^H$ with its image under the Hilbert map $p$.
By Eq.~\eqref{eq:coordinates}, this entails that
\begin{equation}
\label{eq:tangentlimits}
  \widehat{\tau} =
  \lim_{i\rightarrow \infty} T_{\widehat{\Psi} (G(k_i,y_i))} \widehat{S} =  \lim_{i\rightarrow \infty} T_{p(w_i)}p(W_K) \times V^H .
\end{equation}
Since $p(W_K)$ is semialgebraic and an analytic manifold, \cite[Prop.~3, p.~103]{LojESA} by {\L}ojasiewicz entails that
$p(W_K)$ satisfies Whitney's condition B over the origin. This means after possibly passing to subsequences,
that $\ell_{W_K} \subset \tau_{W_K}$, where $\ell_{W_K}$ is the limit line of the secants $\overline{p(w_i),0}$,
and  $\tau_{W_K}$ the limit of the tangent spaces $T_{p(w_i)}p(W_K)$ for $i\rightarrow \infty$.
By Eqs.~\eqref{eq:coordinates} and \eqref{eq:tangentlimits} this entails that
\[
  \widehat{\ell} \subseteq \ell_{W_K}  \times V^H  \subseteq \tau_{W_K} \times V^H = \widehat{\tau} .
\]
This finishes the proof.
\end{proof}

Recall that $\widehat{\rho} : \inertia{M} \to \inertia{X}$ denotes
the quotient map, which is both open and closed by
\cite[Prop.~3.1 (iv) and Prop.~3.6 (i)]{tomdieck}.
Hence, as the sets defining the $\mathcal{S}_{(h,x)}$ consist of entire $G$-orbits,
and as the pieces of $\mathcal{Z}$ consist of connected components of $G$-orbits,
Proposition \ref{prop-StratConnectedFrontier} extends to the local decomposition
in $\inertia{X}$ given by the $\mathcal{R}_{G(h,x)}$.  Therefore,
combining Propositions \ref{prop-StratConnectedManifolds}, \ref{prop-StratConnectedGermsCoincide},
\ref{prop-StratConnectedFrontier}, and \ref{prop-StratConnectedWhitneyB},
we have completed the proof of Theorem \ref{thrm-StratConnected}.

%
%

\section{A De Rham Theorem for the Inertia Space}
\label{sec-deRham}

In this section, we prove a de Rham theorem for the inertia space $\inertia{X}$
analogous to that of \cite{SjamaarDeRham} for singular symplectic reduced spaces.


\subsection{Differential Forms on the Inertia Space}
\label{subsec-DiffForms}

Before constructing differential forms on the inertia space, let us briefly recall from 
\cite[Prop.~1.2.7]{PflaumBook} that a stratification (in the sense of Mather \cite{MatSM}) 
of a locally compact topological space $X$  induces a uniquely determined coarsest decomposition
of $X$ into strata. Applied to our situation, where we consider a compact Lie group $G$ acting 
on a smooth manifold $M$, we thus obtain a coarsest decomposition  $\dec$ of $\inertia{M}$ 
which induces the stratification from Theorem \ref{thrm-StratConnected}. The elements of 
$\dec$ are the strata of $\inertia{M}$. It is easy to see that each stratum from
$\dec$ is $G$-invariant and that the family of quotients 
$\{ G\backslash Z \mid  Z \in \dec \}$ forms a decomposition of $\inertia{X}$
which induces the  natural stratification of the inertia space from Theorem 
\ref{thrm-StratConnected}. Let us introduce some notation: $\iota\co \inertia{M} \to G \times
M$ denotes the natural embedding of $\inertia{M}$ as a subspace and
 $\rho\co G \times M \to G\backslash(G \times M)$ the quotient map.  
Moreover, for each $Z \in \dec$, we denote by $\iota_Z\co Z \to
G \times M$  the inclusion and by $\rho_Z\co Z \to G\backslash Z$ the restricted quotient map.

Let us  construct in the following the sheaf of differential forms on the inertia space. 
Given $k\in \N$ we denote by  $\Omega^k_\inv$  the sheaf of $G$-invariant differential
$k$-forms on $G \times M$ treated as a sheaf on $G\backslash(G \times M)$.  That is, if $U$ 
is an open subset of $G\backslash(G \times M)$, then $\Omega^k_\inv (U)$ consists of 
the differential $k$-forms $\omega \in \Omega^k(\rho^{-1}(U))$ on 
$\rho^{-1}(U) \subseteq G \times M$ such that $L_g^*\omega = \omega$
for all $g \in G$, where $L_g :G \times M \rightarrow G \times M$ denotes the left action by $g$ on $G \times M$.  
Similarly, we let $\Omega_\bas^k$ denote the
subsheaf of $\Omega^k_\inv$ consisting of $G$-basic differential forms
on $G \times M$ or any of the $G$-manifolds $Z \subseteq G\times M$. 
More precisely, $\Omega_\bas^k(U)$ consists of all
$G$-invariant $k$-forms $\omega$ on $\rho^{-1}(U) $ such that 
the interior product $i_{\xi_{G \times M}}\omega $ of $\omega$ with the fundamental vector field 
$\xi_{G \times M}$ vanishes for every $\xi \in \mathfrak{g}$ (cf.~\cite[Sec.~5.3.1]{PflaumBook}).
Now let $W \subseteq \inertia{X}$ be relatively open, and $U\subseteq G\backslash(G \times M)$
open such that $W = U \cap  \inertia{X}$. By a \emph{differential $k$-form} $\widetilde{\omega}$ on
$W$ we now understand  a collection of differential forms $\widetilde{\omega}_Z$ on 
$W\cap(G\backslash Z)$ for $Z \in \dec$ with $W\cap(G\backslash Z)\neq\emptyset$ such that there is 
an $\omega \in \Omega^k_\inv (U)$ with $\rho_{Z}^\ast \widetilde{\omega}_Z =
\iota_{Z}^\ast \omega$ on its domain $\rho^{-1}(W)\cap Z$. We denote the space of differential $k$-forms
on $W$ by $\Omega^k (W)$. One checks immediately that $\Omega^k$ then becomes a sheaf on $\inertia{X}$.  
This sheaf is even fine, since by construction $\Omega^k$ is a $\calC^\infty_X$-module sheaf,
and $\calC^\infty_X$ is fine as the structure sheaf of a differentiable space. 

Note that the form $\omega$ on $U$ which represents the differential form $\widetilde{\omega}$ on
$W$ need not be globally basic.  We let $\Omega_\ibas^k$ denote the subsheaf of $\Omega^k_\inv$ consisting of
$k$-forms $\omega$ such that for every $Z \in \dec$ the pull-back 
$\iota_Z^\ast \omega $ is a basic form on $Z$.  
That is, for each $U \subseteq G\backslash(G \times M)$ open, we define
\[
    \Omega_\ibas^k(U) = \{ \omega \in \Omega^k(\rho^{-1}(U))^G \mid
    i_{\xi_Z} \iota_Z^\ast \omega = 0 \text{ for all $\xi \in \frak{g}$ and 
    $Z \in \dec$}\}.
\]
We refer to sections of $\Omega_\ibas^k$ as \emph{inertia-basic} $k$-forms. Intuitively, these 
correspond to $k$-forms that are basic on each of the strata of $\inertia{M}$.  
A form $\omega \in \Omega_\bas^k(G \backslash ( G \times M) )$ is inertia-basic, 
but an inertia-basic form need not be basic on all of $G \times M$.

By definition, it is clear that we have a surjective linear map
\[
    \Omega_\ibas^k(U)   \longrightarrow   \Omega^k(W)
\]
and that this map has kernel
\[
    \mathcal{I}^k(U)    =   \{ \omega \in \Omega^k (\rho^{-1}(U))^G \mid
    \iota_Z^\ast\omega = 0 \text{ for all $Z \in \dec$} \}.
\]
Hence we obtain isomorphisms
\[
    \Omega^k(W) \cong   \Omega_\ibas^k(U)/\mathcal{I}^k(U).
\]
In particular, when $k = 0$,
\[
    \Omega^0(W) \cong     \Omega_\ibas^0(U)/\mathcal{I}^0(U)=
    \mathcal{C}^\infty(\rho^{-1}(U))^G/\mathcal{I}^0(U),
\]
where $\mathcal{I}^0(U)$ is the ideal of $G$-invariant smooth
functions on $\rho^{-1}(U)$ which vanish on $\inertia{M}$.
By its definition in Section \ref{sec-StrucSheaf} the structure sheaf $\calC^\infty_{\inertia{X}}$  of 
$\inertia{X}$ can be naturally identified with the sheaf $\Omega^0$ on $\inertia{X}$.

Next let us show that the exterior derivative maps inertia-basic forms to inertia-basic forms.
Suppose $\omega$ is an inertia-basic $k$-form on $\rho^{-1}(U)$, i.e.~that $\omega \in \Omega_\ibas^k(U)$.  
By Cartan's Magic Formula, we then conclude for each $Z \in
\dec$ which intersects $\rho^{-1}(U)$ and each $\xi \in \mathfrak{g}$ that
\[
    i_{\xi_{Z}} \iota_Z^\ast d \omega =  \iota_Z^\ast i_{\xi_{G \times M}}d \omega = 
    \iota_Z^\ast(-di_{\xi_{G \times M}}\omega + \mathcal{L}_{\xi_{G \times M}}\omega) = 
    -d\iota_Z^\ast i_{\xi_{G \times M}} \omega = -di_{\xi_Z} \iota_Z^\ast   \omega = 0.
\]
Therefore, $d\omega$ is inertia-basic as well, and we obtain a complex of sheaves
\
\begin{equation}
  \label{eq:resinertia}
   0
    \longrightarrow
    \R_{\inertia{X}}
    \longrightarrow
    \calC^\infty_{\inertia{X}} = \Omega^0
    \stackrel{d}{\longrightarrow}
    \Omega^1
    \stackrel{d}{\longrightarrow}
    \Omega^2
    \stackrel{d}{\longrightarrow}
    \cdots \: ,
\end{equation}
where $\R_{\inertia{X}}$ denotes the sheaf of locally constant $\R$-valued functions on $\inertia{X}$.


\subsection{The Poincar\'{e} Lemma for the Inertia Space}
\label{subsec-deRham}

Let us show that the complex of sheaves \eqref{eq:resinertia} is acyclic, or in other words that a
Poincar\'{e} Lemma holds true for forms on the inertia space. 
So suppose that $\omega$ is a $k$-form on $\rho^{-1}(U)$ for some open $U \subseteq G\backslash(G \times M)$,
and that $d\omega \in \mathcal{I}^{k+1}(U)$. 
Choose a slice $V_{(h,x)}$ at $(h,x) \in \rho^{-1}(U)$ according to Proposition \ref{Prop-LocalInertiaScalarMult}
so that $GV_{(h,x)} \subseteq \rho^{-1}(U)$. By possibly shrinking $V_{(h,x)}$ if necessary,
we may assume by the slice theorem that $Z \cap V_{(h,x)}$ is invariant under the action of $t\in (0,1]$ for every 
$Z \in \dec$.
Let $H = Z_{G_x}(h)$ denote the isotropy group of $(h,x)$. Following \cite[Lemmas 5.2.1 and 5.3.2]{PflaumBook}, we define
\[
   \mathcal{H} \co G \times_H V_{(h,x)} \times [0,1] \longrightarrow G \times_H V_{(h,x)}
\]
by setting
\[
    \mathcal{H}([g,(k,y)], t) = [g, (1 - t)(k,y)].
\]
Then $\mathcal{H}$ is a $G$-invariant retraction of $G \times_H V_{(h,x)}$ onto $G \times_H \{ (h,x)\}$
which restricts to a $G$-invariant retraction of $(G \times_H V_{(h,x)}) \cap \inertia{M}$
onto a single orbit by  Proposition \ref{Prop-LocalInertiaScalarMult}. Let us point out that
by the slice theorem we can naturally identify $G \times_H V_{(h,x)}$ with the set  $GV_{(h,x)} \subseteq \rho^{-1}(U)$.
Next, let $\mathcal{K}: \Omega^k(G \times_H V_{(h,x)} \times [0,1])\to  \Omega^{k-1}(G \times_H V_{(h,x)})$  denote the 
homotopy operator which maps $\omega $ to $\mathcal{K} (\omega)$, where
\[
  \mathcal{K} (\omega) ([g,(k,y)]) =
  \int_0^1 \omega ( [g,(k,y)] , s) \Big( \frac{\partial}{\partial s} , - , \ldots , - \Big) \, ds 
  \text{ for $g\in G$, $(k,y) \in V_{(h,x)}$}.
\]
One checks (see \cite[Lemma 5.2.1]{PflaumBook}), that then
\[
  d \mathcal{K}\mathcal{H}^\ast +  \mathcal{K}\mathcal{H}^\ast d = \mathcal{H}^\ast_1 - \mathcal{H}^\ast_0, 
\]  
where $\mathcal{H}_s = \mathcal{H} ( -,s)$ for $s\in[0,1]$. Hence we obtain for the restriction of $\omega$ to 
$GV_{(h,x)}$
that
\begin{equation}
\label{eq:intclsdform}
  \omega_{|GV_{(h,x)}} - d \mathcal{K}\mathcal{H}^\ast  \omega_{|GV_{(h,x)}} = \mathcal{K}\mathcal{H}^\ast d \omega_{|GV_{(h,x)}}  .
\end{equation}
To prove that the right hand side of this equation lies in $\mathcal{I}^{k}(U')$, where $U':= \rho (V_{(h,x)} )$,
we will show that $\mathcal{K}\mathcal{H}^\ast$ maps $\mathcal{I}^k(U')$ into $\mathcal{I}^{k-1}(U')$.
So suppose that $\eta \in \mathcal{I}^k(U')$ which means 
that $\iota_Z^\ast\eta = 0$ on $\rho^{-1}(U')\cap Z$.  Let $\mathcal{H}_Z$ denote the homotopy
\[
   \mathcal{H}_Z\co G \times_H (Z \cap V_{(h,x)}) \times [0,1] \longrightarrow G \times_H (\overline{Z} \cap V_{(h,x)})
\]
given by restricting $\mathcal{H}$. Similarly, let $\mathcal{K}_Z$ denote the restriction of the operator
$\mathcal{K}$ to $\Omega^k(G\times_H(Z\cap V_{(h,x)}))$.  Then the diagram 
\[
\begin{diagram}[width=2.5cm]
    G \times_H (Z \cap V_{(h,x)}) \times [0,1] &   \rTo^{\iota_Z \times \id_{[0,1]}}
                                    &   (G \times_H V_{(h,x)}) \times [0, 1]          \\
    \dTo<{\mathcal{H}_Z}            &   &   \dTo>{\mathcal{H}}                                    \\
    G \times_H (\overline{Z} \cap V_{(h,x)})           &   \rTo^{\iota_Z}
                                    &       G \times_H V_{(h,x)}                        \\
\end{diagram}
\]
commutes.  Since the operator $\mathcal{K}$ clearly commutes with the restriction to $Z$,
this entails
\[
    \iota_Z^\ast \mathcal{K}\mathcal{H}^\ast \eta
    =  \mathcal{K}_Z\mathcal{H}_Z^\ast \iota_Z^\ast\eta = 0.
\]
Moreover, since  $\mathcal{K}$ and $\mathcal{H}$ commute with the $G$-action,
we obtain for $\xi \in \mathfrak{g}$ 
\[
    i_{\xi_{G V_{(h,x)}}}  \mathcal{K}\mathcal{H}^\ast \eta = \mathcal{K}\mathcal{H}^\ast i_{\xi_{G V_{(h,x)}}} \eta = 0.
\]
It follows that $\mathcal{K}\mathcal{H}^\ast$ maps 
$\mathcal{I}^k(U')$ into $\mathcal{I}^{k-1}(U')$, so that the right hand side of 
Eq.~\eqref{eq:intclsdform} lies in $\mathcal{I}^{k}(U')$, since $d\omega\in \mathcal{I}^{k-1}(U)$ by hypothesis.
But this means that the sheaf complex $\Omega^\bullet$ on 
$\inertia{X}$ is exact, or in other words that the Poincar\'e Lemma for forms on the inertia space holds true.

\begin{theorem}
\label{thrm-deRhamInertia}
The cohomology of the complex $\Omega^\ast(\inertia{X})$ of
differential forms on $\inertia{X}$ naturally coincides with the singular (or \v{C}ech) cohomology
of $\inertia{X}$. Moreover, if $X$ is compact, the cohomology of the de Rham complex $\Omega^\ast(\inertia{X})$
on the inertia space is finite dimensional. 
\end{theorem}
\begin{proof}
By the Poincar\'e Lemma for forms on the inertia space, $\Omega^\bullet$ provides a fine resolution of the sheaf of 
$\R$-valued locally constant functions on $\inertia{X}$. 
Since $\inertia{X}$ is locally compact and locally contractible, the cohomology of the complex $\Omega^\bullet(\inertia{X})$
of global sections then has to coincide naturally with the singular cohomology on $\inertia{X}$.
Since   $\inertia{X}$ is even triangulable, the cohomology of $\Omega^\ast(\inertia{X})$ even coincides with the
\v{C}ech cohomology. 
The triangulability of $\inertia{X}$ also implies that for every open covering of $\inertia{X}$ there exists a 
locally finite subordinate good covering (see \cite[Sec.~7]{PflPosTanGOSPLG}). 
This implies that under the assumption that $X$, hence $\inertia{X}$ is compact, the \v{C}ech cohomology of 
$\inertia{X}$ has to be finite dimensional. This completes the proof.
\end{proof}


%
%




\begin{thebibliography}{10}

\bibitem[{\sc AdGo}]{ademgomez}
\textsc{Adem, A.}, and \textsc{J.M.~G\'{o}mez}:
\emph{Equivariant {K}-theory of compact Lie group actions with maximal rank isotropy},
\texttt{arXiv:1203.4748v1 [math.AT]}

\bibitem[{\sc AdLeRu}]{ademleidaruan}
\textsc{Adem, A., J. Leida}, and \textsc{Y. Ruan},
\emph{Orbifolds and stringy topology},
Cambridge Tracts in Mathematics \textbf{171},
Cambridge University Press, Cambridge, 2007.

\bibitem[{\sc BaCo}]{BaumConnes}
\textsc{Baum, P.}, and \textsc{ A. Connes}:
\emph{Chern character for discrete groups},
A f\'{e}te of topology, 163--232, Academic Press, Boston, MA, 1988.

\bibitem[{\sc BaBrMPh}]{baumbryMacP}
\textsc{Baum, P., J.-L.~Brylinski}, and \textsc{R. MacPherson}:
\emph{Cohomologie \'{e}quivariante d\'{e}localis\'{e}e},
C. R. Acad. Sci. Paris S\`{e}r. I Math. \textbf{300} (1985), 605--608.

\bibitem[{\sc Bie75}]{BierstoneLifting}
\textsc{Bierstone, E.:}
\emph{Lifting isotopies from orbit spaces},
Topology \textbf{14} (1975), 245--252.

\bibitem[{\sc Bie80}]{BierstoneOrbitSpace}
\bysame ,
\emph{The structure of orbit spaces and the singularities of equivariant mappings},
Monograf\'{i}as de Matem\'{a}tica 35. Instituto de Matem\'{a}tica Pura e Aplicada,
Rio de Janeiro, 1980.

\bibitem[{\sc BlGe}]{BlockGetz}
\textsc{Block, J.}, and \textsc{E. Getzler}:
\emph{Equivariant cyclic homology and equivariant differential forms},
Ann. Sci. \'{E}cole Norm. Sup. (4) \textbf{27} (1994),493--527.

\bibitem[{\sc BoCoRo}]{BocCosRoyRAG}
 \textsc{ Bochnak, J., M.~Coste}, and \textsc{M.-F.~Roy}:
 \emph{Real Algebraic Geometry}, Ergebnisse der Mathematik und ihrer
  Grenzgebiete, 3.~Folge, Vol.~\textbf{36}, Springer, Berlin, Heidelberg, New York, 1998.



\bibitem[{\sc Bre}]{BredonBook}
\textsc{Bredon, G.E.}:
Introduction to compact transformation groups,
Pure and Applied Mathematics, Vol.~\textbf{46}. Academic Press, New York-London, 1972.

\bibitem[{\sc BrDi}]{BroeckertomDieck}
{\sc Br\"{o}cker, T.,} and {\sc T. tom Dieck}, \emph{Representations of compact
{L}ie groups}, Graduate Texts in Mathematics, Springer-Verlag, New
York, 1985.

\bibitem[{\sc Bry}]{bryl87}
{\sc Brylinski, J.-L.:}
\emph{Cyclic homology and equivariant theories},
Ann. Inst. Fourier (Grenoble) \textbf{37} (1987), 15--28.






\bibitem[{\sc DeKn}]{DelKneLSS}
 \textsc{H.~Delfs}, and \textsc{M.~Knebusch}:
 \emph{Locally Semialgebraic Spaces},
 Lecture Notes in Mathematics \textbf{1173}, Springer, 1986.

\bibitem[{\sc DuKo}]{duistermaatkolk}
{\sc Duistermaat J.J}, and {\sc J.A.C.~Kolk},
\emph{Lie groups},
Springer-Verlag, Berlin, 2000.


\bibitem[{\sc Far92}]{FarsiKIndex}
{\sc Farsi, C.:}
\emph{${K}$-theoretical index theorems for orbifolds},
Quart.~J.~Math.~Oxford Ser.~(2) \textbf{43} (1992), 183--200.

\bibitem[{\sc Far07}]{FarsiRelIndex}
{\sc Farsi},
\emph{An orbifold relative index theorem},
J.~Geom.~Phys.~\textbf{57} (2007), 1653--1668.

\bibitem[{\sc GoHoKn}]{GoHoKnu07}
{\sc Goldin, R., S.  Holm}, and {\sc A.~Knutson},
\emph{Orbifold cohomology of torus quotients},
Duke Math.~J.~\textbf{139} (2007), 89--139.

\bibitem[{\sc Gor}]{gor}
{\sc Goresky, R.M.:} \textit{Triangulation of Stratified Sets}.
Proc. Amer. Math. Soc.~\textbf{72}, Nr.~1, 193--200 (1980).

\bibitem[{\sc Kaw78}]{kawasaki1}
{\sc Kawasaki, T.:}
\emph{The signature theorem for $V$-manifolds},
Topology \textbf{17} (1978), 75--83.

\bibitem[{\sc Kaw79}]{kawasaki2}
\bysame ,
\emph{The Riemann--Roch theorem for complex $V$-manifolds},
Osaka J. Math.~\textbf{16} (1979), 151--159.

\bibitem[{\sc Kaw84}]{kawasaki3}
\bysame ,
\emph{The index of elliptic operators over $V$-manifolds},
Nagoya Math J.~\textbf{84} (1981), 135--157.

\bibitem[{\sc Kos}]{Koszul}
{\sc Koszul, J.L.:}
\emph{Sur certains groupes de transformation de Lie},
Colloque de G\'{e}om\'{e}trie Diff\'{e}rentielle,
Collogues du CNRS (1953), 137--141.

\bibitem[{\sc Loj}]{LojESA}
  {\sc {\L}ojasiewicz, S}:
  \emph{Ensemble semi-analytique}, Mimeographi\'e, Institute des Hautes
  \'Etudes Scientifique, Bures-sur-Yvette, France, 1965.

\bibitem[{\sc LuUr}]{LupUri}
\textsc{Lupercio, E.,} and \textsc{B. Uribe}:
\emph{Inertia orbifolds, configuration spaces and the ghost loop space},
Q.~J.~Math.~\textbf{55} (2004), no. 2, 185--201.

\bibitem[{\sc Mat70}]{MatNTS}
\textsc{Mather, J.}: \emph{Notes on Topological Stability},
Mimeographed Lecture Notes, Harvard, 1970.

\bibitem[{\sc Mat73}]{MatSM}
\bysame ,
Stratifications and mappings, Dynamical Systems (M. M.
Peixoto, ed.), Academic Press, 1973, pp. 195--232.

\bibitem[{\sc MaSh}]{MatumotoShiota}
{\sc Matumoto, T.}, and {\sc M. Shiota:}
\emph{Proper subanalytic transformation groups and unique triangulation of the orbit spaces},
Transformation groups, Poznan 1985, 290--302, Lecture Notes in Mathematics \textbf{1217}, Springer-Verlag,
Berlin, 1986.

\bibitem[{\sc Mic}]{michor}
{\sc Michor, P.W.:}
\emph{Isometric actions of {L}ie groups and invariants},
Lecture Notes, \texttt{http://www.mat.univie.ac.at/~michor/tgbook.pdf} (1996).

\bibitem[{\sc MoMr}]{moerdijkmrcun}
  {\sc Moerdijk, I.} and {\sc Mr\v{c}un, J.}:
  \textit{Introduction to foliations and Lie groupoids.}
  Cambridge Studies in Advanced Mathematics, \textbf{91}. Cambridge University
  Press, Cambridge, (2003).


\bibitem[{\sc GoSa}]{NGonzalezSanchoBook}
{\sc J.A. Navarro Gonz\'{a}lez}, and {\sc J.B. Sancho de Salas:}
\emph{$C^\infty$-differentiable spaces},
Lecture Notes in Mathematics \textbf{1824}. Springer-Verlag, Berlin, 2003.


\bibitem[{\sc Pfl}]{PflaumBook}
{\sc Pflaum, M.J.:}
\emph{Analytic and geometric study of stratified spaces},
Lecture Notes in Math.~\textbf{1768}, Springer-Verlag, Berlin, 2001.

\bibitem[{\sc PfPoTa07}]{PflaumAlgIndex}
{\sc Pflaum, M.J., H.B.~Posthuma}, and {\sc X.~Tang:},
\emph{An algebraic index theorem for orbifolds},
Adv.~Math.~\textbf{210} (2007), 83--121.

\bibitem[{\sc PfPoTa11}]{PflPosTanGOSPLG}
\bysame ,
\emph{Geometry of orbit spaces of proper Lie groupoids},
\texttt{arXiv:1101.0180v3 [math.DG]}



\bibitem[{\sc Sch}]{Schwarz}
{\sc Schwarz, G.W.:}
\emph{Lifting smooth homotopies of orbit spaces},
Inst. Hautes \'{E}tudes Sci. Publ. Math.~\textbf{51} (1980),  37–135

\bibitem[{\sc Seg}]{Segal}
{\sc Segal, G.B.:}
\emph{The representation ring of a compact Lie group},
Publ. Math. IHES.~\textbf{34} (1968), 113--128

\bibitem[{\sc Sja}]{SjamaarDeRham} {\sc Sjamaar, R.:} \emph{A de
    Rham theorem for symplectic quotients}, Pacific J. Math.~\textbf{220} (2005), 153--166.

\bibitem[{\sc Spa69}]{SpallekDR}
{\sc Spallek, K.:}
Differenzierbare R\"{a}ume,
Math.~Ann.~\textbf{180} (1969), 269--296.

\bibitem[{\sc Spa70}]{SpallekGlaettungDR}
\bysame ,
Gl\"{a}ttung differenzierbarer R\"{a}ume,
Math.~Ann.~\textbf{186} (1970), 233--248.

\bibitem[{\sc Spa71}]{SpallekDiffForms}
\bysame ,
Differential forms on differentiable spaces,
Rend.~Mat.~({6}) \textbf{4} (1971), 231--258.

\bibitem[{\sc Spa72}]{SpallekDiffFormsII}
\bysame ,
Differential forms on differentiable spaces. II,
Rend.~Mat.~({6}) \textbf{5} (1972), 375--389.

\bibitem[{\sc tDie}]{tomdieck}
{\sc tom Dieck, T.:}
\emph{Transformation groups},
de Gruyter Studies in Mathematics 8, Walter de Gruyter, Berlin, 1987.

\bibitem[{\sc Tro}]{Trofimov}
{\sc Trofimov, V.V.:}
\emph{Introduction to geometry of manifolds with symmetry},
Mathematics and its Applications 270, Kluwer Academic Publishers, Dordrecht, The Netherlands, 1994.

\bibitem[{\sc Verg}]{VergneEquivarIndex}
{\sc Vergne, M.:}
\emph{Equivariant index formulas for orbifolds},
Duke Math. J.~\textbf{82} (1996), 637--652.

\bibitem[{\sc Vero}]{verona}
  {\sc Verona, A.:} \textit{Triangulation of Stratified Fibre Bundles}.
  Manuscripta Math.~\textbf{30}, 425--445 (1980).


\end{thebibliography}
\end{document}